\documentclass{amsart}

\usepackage{amssymb}

\newtheorem{theorem}{Theorem}[section]
\newtheorem{lemma}[theorem]{Lemma}

\newtheorem{proposition}[theorem]{Proposition}

\newtheorem{remark}[theorem]{Remark}

\numberwithin{equation}{section}
\newcommand{\Ai}{\text{Ai\,}}

\newcommand{\Tr}{\text{Tr\,}}
\newcommand{\re}{\text{Re\,}}
\newcommand{\im}{\text{Im\,}}

\begin{document}
\setcounter{page}{1}

\thanks{Supported by the G\"oran Gustafsson Foundation (KVA)}

\title[From Gumbel to Tracy-Widom]
{From Gumbel to Tracy-Widom}
\author[K.~Johansson]{Kurt Johansson}

\address{
Department of Mathematics,
Royal Institute of Technology,
S-100 44 Stockholm, Sweden}

\email{kurtj@kth.se}

\begin{abstract}
The Tracy-Widom distribution that has been much studied in recent years can be thought of 
as an extreme value distribution. We discuss interpolation between the classical extreme 
value distribution $\exp(-\exp(-x))$, the Gumbel distribution and the Tracy-Widom
distribution. There is a family of determinantal processes whose edge behaviour
interpolates between a Poisson process with density $\exp(-x)$ and the Airy kernel
point process. This process can be obtained as a scaling limit of a grand canonical
version of a random matrix model introduced by Moshe, Neuberger and Shapiro. We also consider
the deformed GUE ensemble, $M=M_0+\sqrt{2S} V$, with $M_0$ diagobal with independent
elements and $V$ from GUE. Here we do not see a transition from Tracy-Widom to Gumbel,
but rather a transition from Tracy-Widom to Gaussian.

\end{abstract}

\maketitle

\section{Introduction and results}\label{sect1}

\subsection{Introduction}

In the random matrix litterature there has been alot of discussion about the 
transition from Poissonian to random matrix eigenvalue statistics, see for example
\cite{FoPCUE},\cite{FGM-G}, \cite{Gu}, \cite{MCIN}, \cite{Pa}.
One motivation comes from disordered systems, and another from
quantum chaos where Poissonian statistics is expected to describe the eigenvalue 
statistics of classically integrable systems (Berry-Tabor conjecture), and random 
matrix statistics should describe eigenvalue statistics of systems whose classical 
dynamics is fully chaotic (Bohigas-Gianonni-Schmidt conjecture). Hence it has
been natural to look at transitions
between Poissonian and random matrix statistics. In general there could be many
ways to go between different ensembles, but it is nevertheless interesting to find
natural interpolating ensembles and investigate their properties. Mathematically it
is easiest to consider Hermitian (Unitary) ensembles. 
Previous papers on the problem 
have been mainly concerned with 
the transition statistics for eigenvalues in the bulk of the 
spectrum. In the bulk we should see a transition from a Poisson process to a sine kernel
determinantal point process, and for the nearest neigbour spacing statistics we should
see a transition from the exponential distribution to the Gaudin distribution.

In this paper we will discuss the edge behaviour of the eigenvalues. In a finte random
matrix ensemble we look at a scaling limit 
around the largest eigenvalue instead of in the bulk of the 
spectrum. If we take a diagonal matrix with independent Gaussian entries, the largest
eigenvalue will, as the size of the matrix grows, fluctuate according to the Gumbel extreme
value distribution. If we take a full Gaussian matrix from GUE, then the largest eigenvalue
will fluctuate according to the Tracy-Widom distribution. Can we find interesting
distributions that interpolate between Gumbel and Tracy-Widom? Is there a family of
determinantal processes that interpolates? 
Should we typically expect to see a transition from Tracy-Widom to Gumbel.
To shed some light on these question we will
discuss two interpolating random matrix ensembles on Hermitian matrices.

\subsection{The ensembles}

\noindent
(1) {\it Deformed GUE}.
Let $\text{diag\,}(y_1,\dots,y_N)$ denote the diagonal matrix with elements $y_1,\dots,y_N$.
We take $y_1,\dots,y_N$ to be independent Gaussian, $\text{N}\,(0,1/2)$ say, random
variables. Let $V$ be an independent $N\times N$ GUE matrix with density $Z_N^{-1}\exp[-\Tr
V^2]dV$. Consider the random matrix
\begin{equation}\label{1.1}
M=\text{diag\,}(y_1,\dots,y_N)+\sqrt{2S}V,
\end{equation}
where $S\ge 0$ is a parameter. When $S=0$ we have a diagonal matrix with independent
entries and when $S\to\infty$, the matrix $M/\sqrt{S}$ approaches a GUE matrix.

\noindent
(2) {\it MNS-model}.
This model was introduced by Moshe, Neuberger and Shapiro in \cite{MNS}, and we will
call it the {\it MNS-model}. Let $H$ be a Hermitian matrix and $U$ a fixed unitary matrix.
A probability measure on the space of Hermitian matrices is defined by
\begin{equation}
P_{N,U}(H)dH=\frac 1{Z_N}e^{-\Tr H^2}e^{-b\Tr([U,H][U,H]^\ast)}dH,
\notag
\end{equation}
where $b>0$,  $dH$ is Lebesgue measure on the Euclidean space of Hermitian matrices, $[U,H]=
UH-HU$ and the star denotes Hermitian conjugate. The weight is maximal when $[U,H]=0$
so that $U$ and $H$ can be simultaneously diagonalized. The unitary matrix selects
a preferred basis. We can get a unitarily invariant measure by averaging
over the unitary group with respect to the Haar measure, so that we consider a random 
preferred basis. We obtain the probability measure
\begin{equation}\label{1.2}
P_{N}(H)dH=\frac 1{Z_N'}e^{-\Tr H^2}\left(\int_{U(N)}e^{-b\Tr([U,H][U,H]^\ast)}dU\right)dH.
\end{equation}
The integral over the unitary group can be evaluated using the Harish-Chandra or 
Itzykson/Zuber integral and this makes it possible to compute the eigenvalue measure induced
by (\ref{1.2}), \cite{MNS}. This gives
\begin{equation}\label{1.3}
p_N(x)d^Nx=\frac 1{Z_N''}\det\left(e^{-(b+1/2)(x_i^2+x_j^2)+2bx_ix_j}\right)_{1\le i,j\le N}
d^Nx,
\end{equation}
where $x_1,\dots,x_N$ are the eigenvalues of $H$. Actually we will consider a grand 
canonical version of the model, see below, as was also done in \cite{MNS}.
There is a generalization of the MNS-model to Laguerre/Chiral type ensembles that
we will not discuss here, see \cite{G-GV}.

Both of the models above have interpretations in terms of non-interesecting paths.

\noindent
(1) {\it Deformed GUE}.
Consider $N$ standard Brownian motions on the real line, we think of them as particles, 
started at $y_1,\dots,y_N$ at time
0 and conditioned never to intersect. Let $x_1,\dots, x_N$ be the positions of the
particles at time $S$. The probability distribution of $x_1,\dots, x_N$ is the same 
as the eigenvalue distribution of $M$ in (\ref{1.1}), see for example \cite{JoUni}.

\noindent
(2) {\it MNS-model}.
Consider $N$ standard Brownian motions on the real line started at $x_1,\dots, x_N$
at time 0, conditioned to come back to $x_1,\dots, x_N$ at time $t$ and without having had 
any collisions during this time. Put an initial density $\prod_{i=1}^N e^{-x_i^2}$ on
the points $x_1,\dots, x_N$. By a theorem of Karlin and McGregor, \cite{KaMc}, we get
a probability density
\begin{equation}\label{1.4}
\frac 1{\mathcal{Z}_N}\det\left(e^{-\frac 12(x_i^2+x_j^2)-\frac 1{2t}(x_i-x_j)^2}
\right)_{1\le i,j\le N}
d^Nx,
\end{equation}
on the $x_j\,$:s. This is the same as (\ref{1.3}) if we take $b=1/2t$. 
We can think of this as a model of non-intersecting paths on a cylinder.
As stated above we will be interested in the transition at the edge of the spectrum.
The transition in the bulk of the spectrum in the the MNS-model occurs when
$b/N^2\sim c$ or $t\sim 1/2cN^2$, $c>0$ a constant, as $N\to\infty$. (The bulk transition
in deformed GUE occurs when $S\sim C/N^2$.) It is remarked, but not discussed further,
in \cite{MNS} that when the bulk transition occurs, the behaviour at the edge is still like 
that of independent eigenvalues. Below, we will see that there is a transition at the
edge when $b/N^{2/3}\sim c$ as $N\to\infty$. 

\subsection{The Gumbel and Tracy-Widom distributions}

Consider $N$ independent random variables $X_1,\dots,X_N$ with distribution $N(0,1/2)$. Then
it is well known that, \cite{LLR},
\begin{equation}\label{1.5}
\mathbb{P}\left[\frac{\max (X_1,\dots,X_N)-a_N}{b_N}\le x\right]\to F_G(x)\doteq
e^{-e^{-x}}
\end{equation}
as $N\to\infty$, where
\begin{equation}\label{1.6}
a_N=\sqrt{\log N}-\frac{\log(4\pi\log N)}{4\sqrt{\log N}},
\end{equation}
\begin{equation}\label{1.7}
b_N=\frac 1{2\sqrt{\log N}}.
\end{equation}
The distribution function $F_G$ is often called the {\it Gumbel distribution}.
If we think of $X_1,\dots,X_N$ as a point process on the real line with $N$ points and we
take the appropriate scaling limit around the rightmost point we get a Poisson process
on $\mathbb{R}$ with density $e^{-x}$. Its correlation functions are
\begin{equation}\label{1.8}
\rho_k(x_1,\dots,x_k)=\prod_{j=1}^k e^{-x_j},
\end{equation}
$k\ge 1$. 

The {\it Tracy-Widom distribution}, $F_{TW}$, is defined by the Fredholm determinant
\begin{equation}\label{1.9}
F_{TW}(x)=\det(I-K_{\text{Airy}})_{L^2(x,\infty)},
\end{equation}
where
\begin{equation}\label{1.10}
K_{\text{Airy}}(x,y)=\int_0^\infty \Ai(x+\lambda)\Ai(y+\lambda)d\lambda,
\end{equation}
is the {\it Airy kernel}, \cite{TW}. This distribution occurs in several
different places and has 
been much studied in recent years, see \cite{JoICM}, \cite{TW-ICM} for reviews.

If we have a point process on $\mathbb{R}$ then its correlation functions, 
$\rho_k(x_1,\dots,x_k)$, are characterized
by 
\begin{equation}\label{1.10'}
\mathbb{E}\,[\prod_{j}(1+\phi(x_j))]=\sum_{k=0}^\infty\frac 1{k!}\int_{\mathbb{R}^k}
\prod_{j=1}^k\phi(x_j)\rho_{k}(x_1,\dots,x_k) d^kx,
\end{equation}
for any measurable bounded function on $\mathbb{R}$ with compact support. Here the product
in the right hand side is over all particles in the process.

A point process on $\mathbb{R}$ is called determinantal if all its correlation functions, 
$\rho_k(x_1,\dots,x_k)$, $k\ge 1$, exist and are given by
\begin{equation}\label{1.11}
\rho_k(x_1,\dots,x_k)=\det(K(x_i,x_j))_{1\le i,j\le k}
\end{equation}
for some function $K:\mathbb{R}^2\to\mathbb{R}$, the {\it correlation kernel}.
A Poisson process on $\mathbb{R}$ with density $\rho(x)$ can be viewed as a, somewhat 
degenerate, determinantal process with correlation kernel,
\begin{equation}\label{1.11'}
K_{\text{ext}}(x,y)=\begin{cases} 0 & \text{if $x\neq y$} \\
\rho(x) & \text{if $x=y$}\end{cases}
\end{equation}

\subsection{The interpolating process}

Before we discuss the asymptotics of the MNS-model and the deformed GUE model we
will consider a determinantal process which interpolates between the Poisson process
with density $e^{-x}$ and the Airy kernel point process, i.e. the
determinantal process that has kernel (\ref{1.10}). We will see later that this process
can be obtained as a scaling limit of the (grand canonical) MNS-model,
and we will call it the {\it interpolating process}. Also,
we will see that we do not have a transition between the Tracy-Widom 
and the Gumbel distribtions in the deformed GUE ensemble. Rather we will see a transition
from Tracy-Widom to Gaussian. This will be discussed further below.

Define
\begin{equation}\label{1.12}
M_\alpha (x,y)=\int_{-\infty}^\infty\frac{e^{\alpha\lambda}}{e^{\alpha\lambda}+1}
\Ai(x+\lambda)\Ai(y+\lambda)d\lambda.
\end{equation}
That the integral is convergent follows for example from (\ref{2.1}) and the
Cauchy-Schwarz inequality.

\begin{proposition}\label{prop1.1}
The kernel $M_\alpha$ defines a trace class operator in $L^2(a,\infty)$ for any real $a$,
and there is a determinantal process, the interpolating process 
with correlation kernel $M_\alpha$. 
\end{proposition}

That the kernel $M_\alpha$ interpolates between the correlation kernels for the
Poisson process with density $e^{-x}$ and the Airy kernel point process is seen in 
the next theorem. The theorem will be proved in section \ref{sect2}.

\begin{theorem}\label{thm1.2} 
We have the following scaling limits
\begin{equation}\label{1.13}
\lim_{\alpha\to 0+}\frac 1{\alpha} M_{\alpha}(\frac u{\alpha}-\frac 1{2\alpha}\log
(4\pi\alpha^3),\frac v{\alpha}-\frac 1{2\alpha}\log
(4\pi\alpha^3))=K_{\text{ext}}(u,v)
\end{equation}
and
\begin{equation}\label{1.14}
\lim_{\alpha\to\infty}M_{\alpha}(u,v)=K_{\text{Airy}}(u,v).
\end{equation}
\end{theorem}

It is not hard to see that $\int_t^\infty M_\alpha(x,x)dx<\infty$ for any $t$, 
see (\ref{2.2}) below, and hence the 
interpolating process has a last particle almost surely. 
The distribution function $F_\alpha$ 
for this last particle will interpolate between the Gumbel and the Tracy-Widom
distributions.
\begin{theorem}\label{thm1.3} 
The distribution function for the last particle in the interpolating process is
\begin{equation}\label{1.15}
F_\alpha(t)=\sum_{n=0}^\infty\frac{(-1)^n}{n!}\int_{(t,\infty)^n}
\det(M_\alpha(x_i,x_j))_{1\le i,j\le n}d^nx
=\det(I-M_\alpha)_{L^2(t,\infty)}.
\end{equation}
Furthermore
\begin{equation}\label{1.16}
\lim_{\alpha\to 0+}F_\alpha(\frac{\xi}{\alpha}-\frac 1{2\alpha}\log(4\pi\alpha^3))=F_G(\xi)
\end{equation}
and
\begin{equation}\label{1.17}
\lim_{\alpha\to\infty} F_\alpha(\xi)=F_{TW}(\xi).
\end{equation}
\end{theorem}

We postpone the proof to section \ref{sect2}.

There is a different way of obtaining the distribution $F_\alpha$ that is given in
the next proposition, which will be proved in section \ref{sect2}. The construction in
the theorem will
not give us the whole interpolating process though. 

\begin{proposition}\label{prop1.3A}
Let $x_1>x_2>\dots$ be a realization of the Airy kernel point process. Let $y_1,y_2,\dots$
be independent random variables with common distribution function
\begin{equation}\label{1.172}
G_\alpha(x)=\frac{e^{\alpha x}}{1+e^{\alpha x}},
\end{equation}
which are also independent of $\{x_i\}$. Define a new point process by $z_j=x_j+y_j$, 
$j\ge 1$.
Then,
\begin{equation}\label{1.173}
\mathbb{P}[\max_{j\ge 1} z_j\le \xi]=F_\alpha(\xi).
\end{equation}
\end{proposition}

Note that the point process $\{z_j\}$ in the theorem is not the interpolating process, it is 
only the last particle distribution that is the same. Processes with Poissonian edge 
behaviour constructed in a similar way have recently been studied in \cite{RA}.

\subsection{Grand canonical determinantal processes}

The probability measure (\ref{1.4}) on $\mathbb{R}^N$ does not define a finite determinantal
point process on $\mathbb{R}$. To get a determinantal point process we have to consider
a grand canonical ensemble with varying $N$, see for example 
\cite{FoGC}, \cite{JoNIP} for related constructions. Let us first consider a general model 
with the same structure.

Let $X$ be a complete separable metric space with a reference measure $\mu$. Assume that 
$\psi_j$, $j\ge 0$, is an orthonormal family of complex-valued functions in $L^2(X,\mu)$. Also,
let $a_n\ge 0$ be a sequence such that $\sum_{n=0}^\infty a_n<\infty$. Set
\begin{equation}\label{1.18}
\phi(x,y)=\sum_{n=0}^\infty a_n\psi_n(x)\psi_n(y).
\end{equation}
The function $\phi(x,y)$ is well-defined in $L^2(\mu\times\mu)$ and $\phi(x,x)$ is well-defined
in $L^1(\mu)$. We can define a probability measure on $X^n$ by
\begin{equation}\label{1.19}
p_N(x)d^N\mu(x)=\frac 1{Z_N}\det(\phi(x_i,x_j))_{1\le i,j\le N}d^N\mu(x),
\end{equation}
where
\begin{equation}\label{1.20}
Z_N=\int_{X^N}\det(\phi(x_i,x_j))_{1\le i,j\le N}d^N\mu(x).
\end{equation}
Here we assume that $p_N(x)\ge 0$ and $Z_N>0$. We construct a grand canonical point process,
compare \cite{DVJ}, p. 123 , by letting
\begin{equation}\label{1.21}
q_N=\frac{\lambda^N}{N!}\frac{Z_N}{Z(\lambda)}
\end{equation}
be the probability of seeing exactly $N$ particles, and $p_N(x)d^N\mu(x)$ be the probability 
measure for finding particles at $x_1,\dots, x_N$ given that there are exactly $N$ particles. 
Here $Z(\lambda)$ is a normalization constant (grand canonical partition function),
\begin{equation}\label{1.22}
Z(\lambda)=\sum_{N=0}^\infty\frac{\lambda^N}{N!} Z_N,
\end{equation}
where $Z_0=1$. If $g$ is a function in $L^\infty$ with bounded support, then
\begin{equation}\label{1.23}
\mathbb{E}\,[\prod_j (1+g(x_j))]=\sum_{N=0}^\infty\frac{q_N}{Z_N}
\int_{X^N}\prod_{j=1}^N g(x_j)\det(\phi(x_i,x_j))_{1\le i,j\le N}d^N\mu(x),
\end{equation}
where the product in the left hand side is over all particles in the process.
The next theorem, that will be proved in section \ref{sect3}, says that this construction leads 
to a determinantal process.

\begin{theorem}\label{thm1.4}
The grand canonical point process defined above is a determinantal process on $X$ with
correlation kernel
\begin{equation}\label{1.24}
K_\lambda(x,y)=\sum_{n=0}^\infty\frac{\lambda a_n}{1+\lambda a_n}\psi_n(x)\psi_n(y).
\end{equation}
We call this type of process a {\it grand canonical determinantal process}.
\end{theorem}

\subsection{The MNS-model}

In the MNS-model we have $X=\mathbb{R}$ and $\mu$ is the Lebesgue measure, 
in the above construction. We take
\begin{equation}\label{1.25}
\phi(x,y)=\phi_t(x,y)=\frac 1{\sqrt{2\pi t}}e^{-(x^2+y^2)/2-(x-y)^2/2t}
\end{equation}
Then the probability measure (\ref{1.4}) is exactly the measure (\ref{1.19}). That we have an 
expansion of the form (\ref{1.18}) follows from the next lemma, which is just a way of writing
Mehler's formula. We will give the details in the beginning of section \ref{sect4}.

\begin{lemma}\label{lem1.5}
Set $\beta_q=\sqrt{\frac{1+q}{1-q}}$. Then
\begin{align}\label{1.26}
&\frac{\sqrt{q}}{(1-q)\sqrt{\pi}}\exp(-\frac 12(x^2+y^2)-\frac q{(1-q)^2}(x-y)^2)
\notag\\
&=\sum_{n=0}^\infty\beta_qq^{n+1/2}h_n(\beta_q)h_n(\beta_qy)\exp(-\frac{\beta_q^2}2 (x^2+y^2)),
\end{align}
where $h_n(x)$, $n\ge 0$, are the normalized Hermite polynomials.
\end{lemma}
If we make the identification $1/2t=q/(1-q)^2$ and define
\begin{equation}\label{1.27}
\psi_n(x)=\sqrt{\beta_q}h_n(\beta_q x)e^{-\beta_q x^2/2},
\end{equation}
then the $\psi_n$, $n\ge 1$ are orthonormal and $\phi=\phi_t$ can be expanded as in 
(\ref{1.18}) with $a_n=q^{n+1/2}$. Theorem \ref{thm1.4} then gives the next theorem.

\begin{theorem}\label{thm1.6}
The grand canonical MNS-model coming from (\ref{1.3}) or (\ref{1.4}) is a determinantal
point process on $\mathbb{R}$ with correlation kernel
\begin{equation}\label{1.28}
K_\lambda(x,y)=\sum_{n=0}^\infty\frac{\lambda q^{n+1/2}}{1+\lambda q^{n+1/2}}
\psi_n(x)\psi_n(y),
\end{equation}
where the three parameters are related by
\begin{equation}\label{1.29}
b=\frac 1{2t}=\frac q{(1-q)^2}.
\end{equation}
\end{theorem}

Write
\begin{equation}\label{1.29'}
q=e^{-\mu}
\end{equation}
and fix a number $N\ge 0$. Note that $\mu\to\infty$ corresponds to $t\to\infty$ and
$\mu\to 0+$ to $t\to 0+$.
If we choose
\begin{equation}\label{1.30}
\lambda=e^{\mu N}-1,
\end{equation}
then $\int_{\mathbb{R}}K_\lambda(x,x)dx\approx N$, so the expected number of
particles in the process is approximately $N$.

The next proposition shows that the kernel $K_\lambda$ interpolates between a point
process defined by $N$ independent Gaussian random variables and GUE as we should expect.
We postpone the proof to section \ref{sect4}.

\begin{proposition}\label{prop1.7}
If we choose $q$ as in (\ref{1.29'}) and $\lambda$ as in (\ref{1.30}), then
\begin{equation}\label{1.31}
\lim_{\mu\to\infty} K_\lambda(x,y)=\sum_{n=0}^{N-1} h_n(x)h_n(y)e^{-(x^2+y^2)/2}
\doteq K_{GUE(N)}
\end{equation}
uniformly for $x,y$ in a compact set, and
\begin{equation}\label{1.31'}
\lim_{\mu\to 0+} K_\lambda(x,y)=\begin{cases} \frac{N}{\sqrt{\pi}} e^{-x^2}
& \text{if $x\neq y$} \\ 0  & \text{if $x=y$}\end{cases}
\end{equation}
pointwise.
\end{proposition}

As mentioned above the bulk transition occurs when $\mu\sim 1/cN$. This is the limit that was
studied and discussed in \cite{MNS}.

\begin{theorem}\label{thm1.8}
Let $\mu=1/cN$, with $c>0$ fixed, and let $\lambda$ be given by (\ref{1.30}). In this case
$\lambda$ is a constant $\lambda=e^{1/c}-1$.
The following limit holds,
\begin{equation}\label{1.32}
\lim_{N\to\infty}\frac{\pi}{2N\sqrt{c}}K_\lambda\left(\frac{\pi x}{2N\sqrt{c}},
\frac{\pi y}{2N\sqrt{c}}\right)=L_c(x,y)\doteq
\int_0^\infty\frac{\cos\pi(x-y)u}{\lambda^{-1}e^{u^2/c}+1} du
\end{equation}
uniformly for $x,y$ in a compact set.
\end{theorem}

The theorem will be proved in section \ref{sect4}.

Thus in this transition region in the bulk of the point pocess we will have a determinantal 
process with correlation kernel $L_c$. Suitable scaling limits will give the sine kernel as 
$c\to 0+$ and a uniform Poisson process as $c\to\infty$

In \cite{MNS} only the following approximate expression 
\begin{equation}\label{1.32'}
L_c(x,y)\approx\frac{\pi c}2 \frac{\sin \pi(x-y)}{\sinh\pi^2c(x-y)/2}
\end{equation}
is given, valid when $c$ is small. At the end of section \ref{sect4} we will sketch 
an argument leading to this approximate formula without discussing the error.

As briefly mentioned in \cite{MNS}, but not really discussed, when $\mu=1/cN$, the
behaviour at the edge is still like that of independent particles, i.e. we get a
Poisson process with density $e^{-x}$. More precisely we have the following theorem, 
which will be proved in section \ref{sect4}.

\begin{theorem}\label{thm1.9}
Let $\mu=1/cN$, $c>0$ fixed and $\lambda=e^{1/c}-1$ as in the previous theorem. Set
$$
a_N(c)=\sqrt{\log N}-\frac{\log(4\pi\log N/\lambda^2c^2)}{4\sqrt{\log N}}
\notag
$$
and $b_N=(2\sqrt{\log N})^{-1}$ as in (\ref{1.7}). Then,
\begin{equation}\label{1.33}
\lim_{N\to\infty} b_NK_\lambda(a_N(c)+b_N\xi,a_N(c)+b_N\eta)=
K_{ext}(\xi,\eta),
\end{equation}
pointwise.
\end{theorem}

To get an intermediate process at the edge we have to pick a larger $\mu$. 
In fact the intermediate process will be exactly the interpolating process with
kernel $M_\alpha$ discussed above. The next theorem will be proved in section \ref{sect4}.

\begin{theorem}\label{thm1.10}
Choose $\mu=\alpha/N^{1/3}$, $\lambda=e^{\alpha N^{2/3}}-1$, where $\alpha>0$ is fixed.
Then,
\begin{equation}\label{1.34}
\lim_{N\to\infty} \frac{\sqrt{\alpha}}{2N^{1/3}}K_\lambda\left(
N^{1/3}\sqrt{\alpha}+\frac{\sqrt{\alpha}}{2N^{1/3}}\xi,
N^{1/3}\sqrt{\alpha}+\frac{\sqrt{\alpha}}{2N^{1/3}}\eta\right)=M_\alpha(\xi,\eta)
\end{equation}
uniformly for $\xi,\eta$ in a compact set.
\end{theorem}

Hence, in the grand canonical version of the MNS-model, we can see a transition between 
Gumbel statistics and Tracy-Widom statistics for the largest eigenvalue.

\subsection{The deformed GUE model}

We turn now to the deformed GUE model (\ref{1.1}). The bulk transition in this and related models
has been discussed for example in \cite{FoPCUE}, \cite{Gu}, \cite{Pa} 
and we will not discuss it here. It occurs for 
$S\sim c/N^2$, which is the same as for the MNS-model. When we look at the edge, the
behaviour of the deformed GUE will be different than that of the MNS-model. We will not
see a transition between Tracy-Widom and Gumbel. If we choose $S=\alpha^2/N^{2/3}$ we
will see a change at the edge behaviour as we vary $\alpha$, but the transition
will be from Tracy-Widom as $\alpha\to\infty$ to Gaussian as $\alpha\to 0+$. Informally
we can interpret this as follows. The eigenvalue distribution is approximately a semicircle
and with $y_1,\dots,y_N$ fixed we would see Tracy-Widom fluctuations. However, the
fluctuations of $y_1,\dots,y_N$ causes the semicircle to fluctuate, that is the position
of the edge fluctuates like a Gaussian. We can think of the semicircle as fluctuating
basically like $\frac 1N\sum_{i=1}^N y_i$, i.e. like a Gaussian. The effect is that the 
largest eigenvalue will fluctuate like a Tracy-Widom random variable plus an 
independent Gaussian. There is some similarity between this problem and the 
random growth model with random parameters studied in \cite{GTW}.

\begin{theorem}\label{thm1.11}
Let $d\mu(t)$ be a probability measure on $\mathbb{R}$ satisfying $\int td\mu(t)=0$,
$\int t^2d\mu(t)=\sigma^2$ and
$\int|t|^7d\mu(t)<\infty$. Let $y_1,\dots,y_N$ be independent random variables with
distribution $d\mu(t)$ and consider the random $N\times N$ matrix
\begin{equation}\label{1.35}
M=\text{diag\,}(y_1,\dots,y_N)+\sqrt{2S}V,
\end{equation}
where $S=\alpha^2/N^{2/3}$ and 
$V$ is an independent GUE matrix with density $Z_N^{-1}\exp(-\Tr V^2)dV$.
Let $\lambda_{\max}^{(N)}$ be the largest eigenvalue of $M$. There is a number $R(N)\sim 
2\alpha N^{1/6}$,
given by (\ref{5.11}) below, which depends on $d\mu$, $\alpha$ and $N$, so that
\begin{equation}\label{1.36}
\lim_{N\to\infty}\mathbb{P}\left[\frac{\lambda_{\max}^{(N)}-R(N)}{\alpha/\sqrt{N}}\le t
\right]=\mathbb{P}[X+Y\le t],
\end{equation}
where $X$ and $Y$ are independent, $X$ has the Tracy-Widom distribution and $Y$ has
distribution $N(0,\sigma^2/\alpha^2)$.
\end{theorem}

If we want to compare with proposition \ref{prop1.3A} we can let $x_1>x_2>\dots$ be a
realization of the Airy kernel point process and $y$ be an independent random
variable with distribution $N(0,\sigma^2/\alpha^2)$. Set $z_j=x_j+y$. Then $\max_{j\ge 1}
x_j=x_1+y$, will be distributed according to the right hand side of (\ref{1.36}).

\begin{remark}
\rm
Another model for the transition between independent eigenvalues and GUE random matrix 
eigenvalues is a band Hermitian matrix with Gaussian elements. Let $m_{ii}$ be independent 
N(0,1/2), $1\le i\le N$ and $\re m_{ij}$, $\im m_{ij}$, $1\le i<j\le N$, and $j-i<b$,
be independent N(0,1/4), for some given $b$, $1\le b\le N$. Set $m_{ij}=0$ for 
$1\le i<j\le N$, and $j-i\ge b$ and $m_{ij}=m_{ji}$. Then $M=(m_{ij})_{1\le i,j\le N}$
is a diagonal matrix when $b=1$ and a GUE matrix when $b=N$. It is conjectured, see
for example \cite{KK} and references therein, that the local bulk statistics, in the
limit $N\to\infty$, shows a transition from a Poissonian to a determinantal sine-kernel
point process when $b\sim cN^{1/2}$, $0<c<\infty$. When is there a transition at the edge?
Based on the results above one might guess that the edge transition takes place for
a larger $b$. Do we see a transition from Gumbel to Tracy-Widom or is there something else
happening in between? It is not easy to approach these problems. Since we are dealing with
the edge and not the bulk it could be that the method of moments, used with great success
in \cite{SoWi} for Wigner matrices, is useful here also.\it
\end{remark}

\begin{remark}
\rm
The comparison of the Tracy-Widom distribution with the Gumbel distribution suggests
that we are thinking of the Tracy-Widom distribution as a kind of extreme value
distribution. One way to motivate this is as follows. Let $w(i,j)$, $(i,j)\in\mathbb{Z}_+^2$,
be i.i.d. geometric random variables, and let $\pi_k^{(N)}$, $k=1,\dots,\binom{2N}N$,
be all up/right paths from $(1,1)$ to $(N,N)$. Set
\begin{equation}
X_k^{(N)}=\sum_{(i,j)\in\pi_k^{(N)}} w(i,j).
\notag
\end{equation}
For $N$ large each $X_k^{(N)}$ is approximately normal. Clearly, the $X_k^{(N)}$
are not independent. The random variable $G(N,N)=\max_k X_k^{(N)}$ is thus a maximum
over dependent random variables each of which is approximately normal. We know, \cite{JoSh},
that $G(N,N)$, appropriately rescaled converges to the Tracy-Widom distribution, which
thus arises as an extreme value distribution for certain dependent random variables.
We are not aware of any last-passage percolation problems that would interpolate
between Tracy-Widom and Gumbel. 

In measures on partitions both the Gumbel and the Tracy-Widom distribution appear, \cite{VY}.
Are there any natural measures on partitions that interpolate in the way that the MNS-model
does?\it
\end{remark}

\section{The interpolating model}\label{sect2}

In this section we will give the proofs of the results for the interpolating determinantal
process with correlation kernel $M_\alpha$. A basic identity that is useful is
\begin{equation}\label{2.1}
\int_{-\infty}^\infty e^{\alpha t}\Ai(x+t)\Ai(y+t)dt=
\frac 1{\sqrt{4\pi\alpha}}e^{-(x-y)^2/4\alpha-\alpha(x+y)/2+\alpha^3/12}
\end{equation}
for $\alpha>0$ and all $x,y$, see for example \cite{Ok}.

\begin{proof} (Proposition \ref{prop1.1}). We first prove that $m_\alpha$ defined by
(\ref{1.12}) is a trace class operator on $L^2(a,b)$ for $-\infty<a<b\le\infty$. Note that
$M$ is symmetric and
\begin{equation}
\sum_{i,j=1}^n z_i\bar{z}_jM_\alpha(x_i,x_j)=\int_{-\infty}^\infty
\frac{e^{\alpha\lambda}}{e^{\alpha\lambda}+1}\left|\sum_{i=1}^n z_i\Ai(x_i+\lambda)\right|^2
d\lambda,
\notag
\end{equation}
for any complex numbers $z_1,\dots,z_N$ and all $x_i$, so $M(x,y)$ is a Hermitian positive
definite function. Hence by \cite{Si}, it suffices to show that
\begin{equation}\label{2.2}
\int_a^\infty M_{\alpha}(x,x)dx<\infty.
\end{equation}
It then follows that $M_\alpha$ defines a trace class operator on $L^2(a,b)$ with
$\Tr M_\alpha=\int_a^b M_\alpha(x,x)dx$. The inequality (\ref{2.2}) follows from the estimate
\begin{align}\label{2.3}
M_{\alpha}(x,x)&=\int_{-\infty}^\infty\frac{e^{\alpha\lambda}}{e^{\alpha\lambda}+1}
\Ai(x+\lambda)^2d\lambda
\notag\\
&\le\int_{-\infty}^\infty e^{\alpha\lambda}\Ai(x+\lambda)^2d\lambda=
\frac 1{\sqrt{4\pi\alpha}}e^{-\alpha x+\alpha^3/12},
\end{align}
by (\ref{2.1}).

If we can show that $0\le M_\alpha\le I$, it follows that there is a determinantal
process with correlation kernel $M_\alpha$, see \cite{SoDet}. Let $f$ be a continuous function
on the real line with compact support. Then,
\begin{equation}
\int_{-\infty}^\infty\int_{-\infty}^\infty M_\alpha(x,y)f(x)f(y)dxdy=
\int_{-\infty}^\infty\frac{e^{\alpha\lambda}}{e^{\alpha\lambda}+1}
\left(\int_{-\infty}^\infty\Ai(x+\lambda)f(x)dx\right)^2d\lambda
\notag
\end{equation}
by Fubini's theorem and hence $M_\alpha\ge 0$. Fix $\epsilon>0$, $0<\epsilon<\alpha$, and 
note that
\begin{equation}\label{2.3'}
\frac{e^{\alpha\lambda}}{e^{\alpha\lambda}+1}\le e^{\epsilon\lambda}
\end{equation}
for all real $\lambda$. Thus
\begin{align}
&\int_{-\infty}^\infty\frac{e^{\alpha\lambda}}{e^{\alpha\lambda}+1}
\left(\int_{-\infty}^\infty\Ai(x+\lambda)f(x)dx\right)^2d\lambda
\notag\\
&\le
\int_{-\infty}^\infty\int_{-\infty}^\infty 
\left(\int_{-\infty}^\infty e^{\epsilon\lambda}\Ai(x+\lambda)\Ai(y+\lambda)\right)f(x)f(y)dxdy
\notag\\
&\int_{-\infty}^\infty\int_{-\infty}^\infty 
\frac 1{\sqrt{4\pi\epsilon}}e^{-(x-y)^2/4\epsilon-\epsilon(x+y)+\epsilon^3/12}f(x)f(y)dxdy.
\notag
\end{align}
Since $f$ is continuous and has compact support this last integral $\to||f||_2^2$ as
$\epsilon\to 0+$. Since
$\epsilon>0$ can be taken arbitrarily small we obtain $M_\alpha\le I$.
\end{proof}

Next we turn to the scaling limits of the kernel $M_\alpha$. 

\begin{proof} (Theorem \ref{thm1.2}).
Set $f(\alpha)=\frac 1{2\alpha}\log(4\pi\alpha^3)$. Then
\begin{align}\label{2.4}
M_\alpha^\ast(u,v)&\doteq\frac 1{\alpha}M_\alpha(\frac u{\alpha}-f(\alpha),
\frac v{\alpha}-f(\alpha))
\notag\\
&=\frac 1{\alpha}\int_{-\infty}^\infty\frac{e^{\alpha(t+f(\alpha))}}
{e^{\alpha(t+f(\alpha))}+1}\Ai(t+u/\alpha)\Ai(t+v/\alpha)dt
\notag\\
&=\sqrt{4\pi\alpha}\int_{-\infty}^\infty\frac{e^{\alpha t}}{\sqrt{4\pi\alpha^3}e^{\alpha t}
+1}\Ai(t+u/\alpha)\Ai(t+v/\alpha)dt.
\end{align}
Using the identity (\ref{2.1}) this can be written as $M_\alpha^\ast(u,v)=A-B$,where
\begin{align}
A&=e^{(u-v)^2/4\alpha^3-(u+v)/2+\alpha^3/12}
\notag\\
B&=4\pi\alpha^2\int_{-\infty}^\infty\frac{e^{2\alpha t}}{\sqrt{4\pi\alpha^3}e^{\alpha t}
+1}\Ai(t+u/\alpha)\Ai(t+v/\alpha)dt.
\notag
\end{align}
It is clear that $A\to 0$ as $\alpha\to 0+$ if $u\neq v$ and $A\to e^{-u}$ as $\alpha\to
0+$ if $u=v$. Hence, we have to show that $B\to 0$ as $\alpha\to 0+$. We can assume that
$u\ge v$ without loss og generality and write
\begin{equation}
B=B_1+B_2+B_3=
\left(\int_{-v}^\infty+\int_{-u}^{-v}+\int_{-\infty}^{-u}\right)
\frac{4\pi\alpha e^{2y}}{\sqrt{4\pi\alpha^3}e^y+1}\Ai(\frac{y+u}{\alpha})\Ai(\frac{y+v}{\alpha})dy.
\notag
\end{equation}
Now $x\to\Ai(x)$ is a bounded function and we have the estimates

\begin{equation}\label{2.5}
|\Ai(x)|\le\frac {C}{|x|^{1/4}}
\end{equation}
for $x<0$, and
\begin{equation}\label{2.6}
|\Ai(x)|\le\frac {Ce^{-2x^{3/2}/3}}{x^{1/4}}
\end{equation}
for $x>0$, where $C$ is a numerical constant.
It follows from these estimates that
\begin{equation}
|B_1|\le C\alpha^{3/2}e^{-2v}\int_0^\infty \frac{e^{2y}}{\sqrt{y}}e^{-4y^{3/2}/
3\alpha^{3/2}}dy
\notag
\end{equation}
Clearly, $B_1\to 0$ as $\alpha\to 0+$. Similarly,
\begin{equation}
|B_2|\le C\alpha^{3/2}\int_{-u}^{-v}\frac{e^{2y}}{|y+u|^{1/4}|y+v|^{1/4}}e^{- 
(y+u)^{3/2}/3\alpha^{3/2}}dy
\notag
\end{equation}
and hence $B_2\to 0$ as $\alpha\to 0+$. Finally,
\begin{equation}
|B_3|\le C\alpha^{3/2}e^{-2u}\int_{-\infty}^0\frac{e^{2y}}{\sqrt{|y|}}dy,
\notag
\end{equation}
which goes to 0 as $\alpha\to 0+$. This proves (i) in the theorem.

The fact that $M_\alpha(u,v)\to K_{Airy}(u,v)$ as $\alpha\to\infty$ follows from 
the estimates (\ref{2.5}), (\ref{2.6}) and the dominated convergence theorem.
\end{proof}

It follows from the estimate (\ref{2.2}) that the interpolating process has a last 
particle almost surely. Its distribution function is given by theorem \ref{thm1.3} 
which we now prove.

\begin{proof} (Theorem \ref{thm1.3}). 
It follows from Hadamard's inequality and (\ref{2.5}) that
\begin{equation}
\sum_{n=0}^\infty\frac 1{n!}\int_{(t,\infty)^n}\det(M_\alpha(x_i,x_j))d^nx
\le\sum_{n=0}^\infty\frac 1{n!}\left(\int_t^\infty M_\alpha(x,x)dx\right)^n<\infty.
\notag
\end{equation}
Hence, the first equality in (\ref{1.15}) holds, see for example \cite{JoHou}. The 
second inequality follows since $M_\alpha$ is trace class on $L^2(t,\infty)$ for any $t$
and $\Tr M_\alpha=\int_t^\infty M_\alpha (t,t)dt$, \cite{GGK}.

Next, we turn to the proof of (\ref{1.16}). Write 
$f(\alpha)=\frac 1{2\alpha}\log(4\pi\alpha^3)$. Then
\begin{equation}\label{2.7}
F_\alpha(\xi/\alpha-f(\alpha))=
\sum_{n=0}^\infty\frac {(-1)^n}{n!}\int_{(\xi,\infty)^n}\det(M_\alpha^\ast(x_i,x_j))d^nx,
\end{equation}
where $M_\alpha$ is as in (\ref{2.4}). It follows from the estimate (\ref{2.3}) that
\begin{equation}
M_\alpha^\ast(x,x)\le\frac 1{\sqrt{4\pi\alpha^3}}e^{-x+\alpha f(\alpha)+\alpha^3/12}
=e^{-x+\alpha^3/12}.
\notag
\end{equation}
Hence, by Hadamard's inequality,
\begin{equation}
\det(M_\alpha^\ast(x_i,x_j))_{1\le i,j\le n}\le e^{n\alpha^3/12}e^{-\sum_{j=1}^n x_j},
\notag
\end{equation}
and it follows from (\ref{1.13}), (\ref{2.7}) and the dominated convergence theorem that
\begin{align}
&\lim_{\alpha\to 0+} F_\alpha(\xi/\alpha-f(\alpha))=
\sum_{n=0}^\infty\frac {(-1)^n}{n!}\int_{(\xi,\infty)^n}\prod_{i=1}^n e^{-x_i}d^nx
\notag\\
&=\sum_{n=0}^\infty\frac {(-1)^n}{n!}e^{-n\xi}=F_G(\xi),
\notag
\end{align}
which proves (\ref{1.16}).

If we use the estimate (\ref{2.3'}) with $\epsilon=1$ and (\ref{2.1}) 
we see that for $\alpha\ge 1$
\begin{equation}\label{2.8}
M_\alpha(x,x)\le e^{-x}.
\end{equation}
Hence, for $\alpha$ large, we have
\begin{equation}
\det(M_\alpha(x_i,x_j))_{1\le i,j\le n}\le e^{-\sum_{j=1}^n x_j}.
\notag
\end{equation}
Consequently, we can use (\ref{1.14}), (\ref{1.15}) and the dominated convergence theorem
to conclude that (\ref{1.17}) holds.
\end{proof}

We will now prove proposition \ref{prop1.3A} which gives an alternative representation
of the $F_\alpha$\,-distribution.

\begin{proof}
(Proposition \ref{prop1.3A}). We have
\begin{align}
\mathbb{P}[\max_{j\ge 1} z_j\le\xi]&=\mathbb{E}\,[\prod_{j=1}^\infty 
(1-\chi_{(\xi,\infty)}(z_j))]=
\mathbb{E}\,[\prod_{j=1}^\infty (1-\chi_{(\xi,\infty)}(x_j+y_j))]
\notag\\
&=\mathbb{E}_x\left[\prod_{j=1}^\infty\left(\int_{-\infty}^\infty (1-\chi_{(\xi,\infty)}(x_j
-y))dG_\alpha(y)\right)\right],
\notag
\end{align}
where $\mathbb{E}_x$ denotes expectation with respect to the Airy kernel point process.
Here we have used the fact that the $y_j$:s are independent with distribution $G_\alpha$
and that they are independent of the Airy kernel point process. 
The last equality then follows from Fubini's theorem. We have also used the fact that the 
$G_\alpha$-distribution is symmetric to replace $y$ with $-y$. 
Note that $\chi_{(\xi,\infty)}(x_j-y)=0$ if and only if $y\ge x_j-\xi$ and thus the last 
expression can be written
\begin{align}
&\mathbb{E}_x\left[\prod_{j=1}^\infty(1-G_\alpha(x_j-\xi))\right]
\notag\\
&=\sum_{k=0}^\infty\frac {(-1)^k}{k!}\int_{\mathbb{R}^k}\prod_{i=1}^kG_\alpha(x_i-\xi)\det
(\int_0^\infty\Ai(x_i+t)\Ai(x_j+t)dt)d^kx
\notag
\end{align}
since we have the Airy kernel point process with correlation kernel (\ref{1.10}).
We now make the shift $x_j\to x_j+\xi$ and manipulate the expressions as follows
\begin{align}
&=\sum_{k=0}^\infty\frac {(-1)^k}{k!}\int_{\mathbb{R}^k}\prod_{i=1}^kG_\alpha(x_i)\det
(\int_\xi^\infty\Ai(x_i+t)\Ai(x_j+t)dt)d^kx
\notag\\
&=\sum_{k=0}^\infty\frac {(-1)^k}{k!}\int_{\mathbb{R}^k}d^kx\prod_{i=1}^kG_\alpha(x_i)
\int_{(\xi,\infty)^k}d^kt\det(\Ai(x_i+t_i)\Ai(x_j+t_i))
\notag\\
&=\sum_{k=0}^\infty\frac {(-1)^k}{k!}\int_{(\xi,\infty)^k}d^kt
\int_{\mathbb{R}^k}d^kx\left(\prod_{i=1}^kG_\alpha(x_i)\Ai(x_i+t_i)\right)
\det(\Ai(x_i+t_j)),
\notag
\end{align}
where we have used the fact that the determinant is unchanged under transposition.
Now, this last expression can be written
\begin{align}
&\sum_{k=0}^\infty\frac {(-1)^k}{k!}\int_{(\xi,\infty)^k}
\det\left(\int_{-\infty}^\infty G_\alpha(x)\Ai(x+t_i)\Ai(x+t_j)dx\right) d^kt
\notag\\
&=\det(I-M_\alpha)_{L^2(\xi,\infty)}=F_\alpha(\xi).
\notag
\end{align}
\end{proof}

\section{The grand canonical point process}\label{sect3}

In this section we will show that the grand canonical point process defined in section
\ref{sect1} using (\ref{1.19}) and (\ref{1.21}) is a determinantal process with correlation 
kernel given by (\ref{1.24}). The proof is based on the identity (\ref{1.23}).

\begin{proof} (Theorem \ref{thm1.4}). We want to prove that
\begin{align}\label{3.1}
\mathbb{E}\,[\prod_{j}(1+g(x_j))]&=
\frac 1{Z(\lambda)}\sum_{N=0}^\infty\frac {\lambda^N}{N!}\int_{X^N}
\prod_{j=1}^N g(x_j)\det(\phi(x_i,x_j))_{1\le i,j\le N}d^N\mu(x)
\notag\\
&=\sum_{N=0}^\infty\frac {1}{N!}\int_{X^N}
\prod_{j=1}^N g(x_j)\det(K_\lambda(x_i,x_j))_{1\le i,j\le N}d^N\mu(x).
\end{align}
The first equality is just (\ref{1.21}) and (\ref{1.23}). The identity (\ref{3.1}) implies 
the theorem, see e.g. \cite{JoHou}.
To prove (\ref{3.1}) we will use some facts on von Koch determinants, see e.g. \cite{GGK}. Let 
$(a_{ij})_{i,j=0}^\infty$ be an infinite matrix and assume that 
\begin{equation}\label{0.0}
\sum_{i=0}^\infty |a_{ii}|<\infty\quad,\quad
\sum_{i,j=0}^\infty |a_{ij}|^2<\infty.
\end{equation}
Then,
\begin{equation}\label{0.1}
\det(I+A)=\sum_{n=0}^\infty\frac 1{n!}\sum_{m\in\mathbb{N}^n}\det(a_{m_im_j})_{i,j=1}^n
\end{equation}
is well-defined. Furthermore, if we have two such matrices $A$ and $B$, then
\begin{equation}\label{0.2}
\det(I+A)\det(I+B)=\det(I+A+B+AB).
\end{equation}
Inserting (\ref{1.18}) into the left hand side of (\ref{3.1}) we get
\begin{align}
&\frac 1{Z(\lambda)}\sum_{N=0}^\infty\frac{\lambda^N}{N!}
\int_{X^N}\prod_{j=1}^N(1+g(x_j))\frac 1{N!}\sum_{m\in\mathbb{N}^N}
\prod_{j=1}^N a_{m_j} \det(\psi_{m_i}(x_j))\det(\psi_{m_i}(x_j))d^Nx
\notag\\
&=\frac 1{Z(\lambda)}\sum_{N=0}^\infty\frac{\lambda^N}{N!}
\sum_{m\in\mathbb{N}^N}\prod_{j=1}^N a_{m_j}\det\left(\int_{X}(1+g(x))\psi_{m_i}(x)
\psi_{m_j}(x)dx\right)
\notag\\
&=\frac 1{Z(\lambda)}\det(I+D_g),
\notag
\end{align}
where
\begin{equation}
D_g(i,j)=\lambda a_i^{1/2}\int_{X}(1+g(x))\psi_i(x_)\psi_j(x)d\mu(x) a_j^{1/2}.
\notag
\end{equation}
Here we have used the determinatal identity
\begin{equation}\label{ex}
\det\left(\int_X\phi_i(x)\psi_j(x)d\mu(x)\right)
=\frac 1{N!}\int_{X^N}\det(\phi_i(x_j))
\det(\psi_i(x_j))d^N\mu(x),
\end{equation}
where all the determinants are of size $N\times N$, see e.g. \cite{JoHou}.
Clearly,
$D_0(i,j)=\lambda a_i\delta_{ij}$,
and hence
$Z(\lambda)=\det(I+A)$,
where
\begin{equation}
A(i,j)=\lambda a_i\delta_{ij}.
\notag
\end{equation}
Set
\begin{equation}
B_g(i,j)=\frac{\lambda a_i^{1/2}}{1+\lambda a_i}\int_{X}g(x)\psi_i(x)\psi_j(x)d\mu(x) 
a_j^{1/2}.
\notag
\end{equation}
Note that since $g$ is bounded we have $|D_g(i,j)|\le Ca_i^{1/2}a_j^{1/2}$,
$|B_g(i,j)|\le Ca_i^{1/2}a_j^{1/2}$, so the conditions (\ref{0.0}) are satisfied.
Note also that,
\begin{equation}
\delta_{ij}+\lambda a_i\delta_{ij}+B_g(i,j)+(AB_g)_{ij}=
\delta_{ij}+D_g(i,j).
\notag
\end{equation}
Hence, by (\ref{0.2})
\begin{align}
&\frac 1{Z(\lambda)}\det(I+D_g)=\frac 1{\det(I+A)}\det(I+A)\det(I+B_g)
\notag\\
&=\det(I+B_g)=\sum_{N=0}^\infty\frac 1{N!}\sum_{m\in\mathbb{N}^N}\det
(B_g(m_i,m_j))_{1\le i,j\le N}
\notag\\
&=\sum_{N=0}^\infty\frac 1{N!}\sum_{m\in\mathbb{N}^N}
\prod_{i=1}^N\frac{\lambda a_{m_i}}{1+\lambda a_{m_i}}\frac 1{N!}\int_{X^N}
\det(\psi_{m_i}(x_j))\det(\psi_{m_i}(x_j))\prod_{j=1}^N g(x_j)d^N\mu(x)
\notag\\
&=\sum_{N=0}^\infty\frac 1{N!}\int_{X^N}\prod_{j=1}^N g(x_j)
\det\left(\sum_{m=0}^\infty\frac{\lambda a_m}{1+\lambda a_m}\psi_m(x_i)\psi_m(x_j)
\right)d^N\mu(x),
\notag
\end{align}
which is the right hand side of (\ref{3.1}). Here we have used the identity (\ref{ex}) again.
\end{proof}

\section{The MNS-model}\label{sect4}

In this section we will give the proofs for the results on the MNS-model stated in
section \ref{sect1}. First we must prove lemma \ref{lem1.5} which makes it possible to
use the formalism for a grand canonical determinantal process and obtain theorem 1.6,
which is the starting point for the asymptotic analysis.

\begin{proof}
(Lemma \ref{lem1.5}).
We will use Mehler's formula,
\begin{equation}\label{4.1}
\sum_{n=0}^\infty\frac{H_n(x)H_n(y)}{2^nn!}q^n=\frac 1{\sqrt{1-q^2}}\exp(
-\frac{q^2}{1-q^2}(x^2+y^2)+\frac{2q}{1-q^2}xy),
\end{equation}
where $0<q<1$ and $H_n$ are the standard Hermite polynomials. If we use
instead the normalized Hermite polynomials,
\begin{equation}\label{4.2}
h_n(x)=\frac 1{\pi^{1/4}\sqrt{2^nn!}}H_n(x)
\end{equation}
and rewrite the exponent we obtain
\begin{equation}\label{4.3}
\sum_{n=0}^\infty q^n h_n(x)h_n(y)e^{-(x^2+y^2)/2}=\frac 1{\pi\sqrt{1-q^2}}\exp(
-\frac{1-q}{2(1+q)}(x^2+y^2)+\frac{q}{1-q^2}(x-y)^2).
\end{equation}
The change of variables $x\to\beta_q x$, $y\to\beta_q y$ and multiplication by 
$\beta_qq^{1/2}$ now gives (\ref{1.26}).
\end{proof}

The choice of the parameter $\lambda$ in (\ref{1.30}) with $q$ given by (\ref{1.29'}) is 
motivated by the fact that the expected number of particles becomes

\begin{align}
\int_{\mathbb{R}}K_\lambda(x,x)dx&=\sum_{n=0}^\infty\frac{\lambda q^{n+1/2}}
{1+\lambda q^{n+1/2}}\approx\int_0^\infty\frac{\lambda e^{-\mu x}}{1+\lambda e^{-\mu x}}dx
\notag\\
&=\frac 1{\mu}\log(1+\lambda)=N.
\notag
\end{align}

We turn to the proof of proposition \ref{prop1.7}

\begin{proof} (Proposition \ref{prop1.7}).
With $q$ as in (\ref{1.29'}) and $\lambda$ as in (\ref{1.30}) we have
\begin{equation}
K_\lambda(x,y)=\sum_{n=0}^\infty\frac 1{1+(1-e^{-\mu N})^{-1}e^{(n+1/2-N)\mu}}
\psi_n(x)\psi_n(y),
\notag
\end{equation}
where $\psi_n$ is given by (\ref{1.27}). We split this into two sums, one from
$n=0$ to $N-1$, called $\Sigma_1$, and one from $N$ to infinity, called $\Sigma_2$.
Since $\beta_q\to 1$ as $\mu\to\infty$ we see that $\Sigma_1$ converges to the right
hand side of (\ref{1.31}). We have to prove that $\Sigma_2\to 0$ as $\mu\to\infty$.
A useful bound is
\begin{equation}\label{4.3'}
|h_n(x)e^{-x^2/2}|\le\frac C{n^{1/12}},
\end{equation}
for all $x$, see \cite{Kr}. Hence,
\begin{equation}
\Sigma_2\le C\sum_{n=N}^\infty e^{(N-n-1/2)\mu} \frac 1{n^{1/6}},
\notag
\end{equation}
which goes to $0$ as $\mu\to\infty$. This proves (\ref{1.31}). 

To prove (\ref{1.31'}) consider first the case $x=y$. We have
\begin{equation}
1\le 1+(e^{\mu N}-1)e^{-(n+1/2)\mu}\le e^{\mu N}.
\notag
\end{equation}
Since
\begin{equation}
K_\lambda(x,y)=\sum_{n=0}^\infty\frac{(e^{\mu N}-1)e^{-(n+1/2)\mu}}
{1+(e^{\mu N}-1)e^{-(n+1/2)\mu}}\psi_n(x)\psi_n(y)
\notag
\end{equation}
we obtain
\begin{align}
&\frac{1-e^{-\mu N}}{\mu N}\mu N\sum_{n=0}^\infty e^{-(n+1/2)\mu}
\psi_n(x)^2\le K_\lambda(x,x)
\notag\\
&\frac{e^{\mu N}-1}{\mu N}\mu N\sum_{n=0}^\infty e^{-(n+1/2)\mu}
\psi_n(x)^2\le K_\lambda(x,x).
\notag
\end{align}
The formula (\ref{1.26}) gives
\begin{equation}
\frac{1-e^{-\mu N}}{\mu N}\frac{\mu N}{(1-e^{-\mu})\sqrt{\pi}}
e^{-\mu/2}e^{-x^2}\le K_\lambda(x,x)
\le\frac{e^{\mu N}-1}{\mu N}\frac{\mu N}{(1-e^{-\mu})\sqrt{\pi}}
e^{-\mu/2}e^{-x^2}.
\notag
\end{equation}
By letting $\mu\to 0+$ we get the first part of (\ref{1.31'}).

Consider now the second case, $x\neq y$. Write
\begin{align}
K_\lambda(x,y)&=(e^{\mu N}-1)\sum_{n=0}^\infty e^{-(n+1/2)\mu}\psi_n(x)\psi_n(y)
\notag\\
&+(e^{\mu N}-1)\sum_{n=0}^\infty \left[\frac 1{1+(e^{\mu N}-1)e^{-(n+1/2)\mu}}-1\right]
e^{-(n+1/2)\mu}\psi_n(x)\psi_n(y)
\notag\\
&\doteq S_1+S_2.
\notag
\end{align}
By (\ref{1.26}),
\begin{equation}
S_1=\frac{e^{\mu N}-1}{1-e^{-\mu}} e^{-\mu/2}e^{-(x^2+y^2)/2-e^{-\mu}(1-e^{-\mu})^{-2}
(x-y)^2},
\notag\end{equation}
which $\to 0$ as $\mu\to 0+$. Furthermore,
\begin{align}
|S_2|&\le (e^{\mu N}-1)^2\sum_{n=0}^\infty e^{-(n+1/2)\mu}|\psi_n(x)\psi_n(y)|
\notag\\
&\le \frac 1{2\sqrt{\pi}} (e^{\mu N}-1)^2\frac {e^{-\mu/2}}{1-e^{-\mu}}(e^{-x^2}+e^{-y^2})
\notag
\end{align}
by the Cauchy-Schwarz' inequality and (\ref{1.26}). We see that the last expression
$\to 0$ as $\mu\to 0+$.
\end{proof}

We turn now to the proof of the two theorems that concern 
the asymptotic behaviour of the kernel $K_\lambda$ in the
regime where we have a transition in the bulk.

\begin{proof} (Theorem \ref{thm1.8}).
We will use the following asymptotic formula for the Hermite polynomials,
\cite{DKMVZ}, valid for $-1+\delta\le x\le 1-\delta$, $\delta>0$ fixed,
\begin{equation}\label{4.4}
h_n(\sqrt{2n}x)e^{-nx^2}=\frac{2^{1/4}}{n^{1/4}\sqrt{\pi}}\frac 1{(1-x^2)^{1/4}}
(\cos[2nF(x)-\frac 12\arcsin x]+O(\frac 1n)),
\end{equation}
where
\begin{equation}
F(x)=\int_x^1\sqrt{1-y^2}dy=\frac12(\arccos x-x\sqrt{1-x^2}).
\notag
\end{equation}
Set $A_N=\lambda^{-1}e^{1/2cN}$. Note that $A_N=\lambda^{-1}+O(1/N)$ and
$\beta_q=\sqrt{2cN}+O(1/N^{3/2})$ as $N\to\infty$.
Write $f_n(x)=h_n(x)e^{-x^2/2}$
Using the asymptotic formula (\ref{4.4}) we obtain
\begin{align}
&\frac{\pi}{2N\sqrt{c}}K_\lambda(\frac{\pi x}{2N\sqrt{c}},\frac{\pi y}{2N\sqrt{c}})=
\frac{\pi\beta_q}{N\sqrt{2c}}\sum_{n=0}^\infty\frac 1{A_Ne^{n/cN}+1}
f_n(\frac{\pi \beta_qx}{2N\sqrt{c}})f_n(\frac{\pi \beta_qy}{2N\sqrt{c}})
\notag\\
&=\frac{\beta_q}{N\sqrt{2c}}\sum_{n=1}^\infty\frac 1{A_Ne^{n/cN}+1}
\left(1-\left(\frac{\pi \beta_qx}{2N\sqrt{2cn}}\right)^2\right)^{-1/4}
\left(1-\left(\frac{\pi \beta_qy}{2N\sqrt{2cn}}\right)^2\right)^{-1/4}
\notag\\
&\times\frac 1{n^{1/2}}
\cos\left[2nF\left(\frac{\pi \beta_qx}{2N\sqrt{2cn}}\right)-\frac 12\arcsin\left(
\frac{\pi \beta_qx}{2N\sqrt{2cn}}\right)\right]
\notag\\
&\times\cos\left[2nF\left(\frac{\pi \beta_qy}{2N\sqrt{2cn}}\right)-\frac 12\arcsin\left(
\frac{\pi \beta_qy}{2N\sqrt{2cn}}\right)\right]+o(1)
\notag\\
&=\frac 1{\sqrt{N}}\sum_{n=1}^\infty\frac 1{\lambda^{-1}e^{n/cN}+1}
\frac 1{n^{1/2}}
\cos\left[2nF\left(\frac{\pi \beta_qx}{2N\sqrt{2cn}}\right)-\frac 12\arcsin\left(
\frac{\pi \beta_qx}{2N\sqrt{2cn}}\right)\right]
\notag\\
&\times\cos\left[2nF\left(\frac{\pi \beta_qy}{2N\sqrt{2cn}}\right)-\frac 12\arcsin\left(
\frac{\pi \beta_qy}{2N\sqrt{2cn}}\right)\right]+o(1).
\notag
\end{align}
Note that the sum
\begin{equation}
\frac 1{\sqrt{N}}\sum_{n=1}^\infty\frac 1{\lambda^{-1}e^{n/cN}+1}
\frac 1{n^{1/2}}
\notag
\end{equation}
is bounded in $N$. We have
\begin{equation}
\frac{\pi \beta_qx}{2N\sqrt{2cn}}=\frac{\pi x}{2\sqrt{nN}}+O(\frac 1{\sqrt{n}N^{5/2}})
\notag
\end{equation}
and
\begin{equation}
2nF\left(\frac{\pi x}{2\sqrt{nN}}\right)=\frac{\pi n}2-\pi x\sqrt{\frac nN}+
O(\frac 1{\sqrt{n}N^{5/2}})
\notag
\end{equation}
as $N\to\infty$, $x$ in a compact set. Hence,
\begin{align}
&\frac{\pi}{2N\sqrt{c}}K_\lambda(\frac{\pi x}{2N\sqrt{c}},\frac{\pi y}{2N\sqrt{c}})=
\frac 1{2N}\sum_{n=1}^\infty\frac{(-1)^n}{\lambda^{-1}e^{n/cN}+1}
\left(\frac nN\right)^{-1/2}\cos(\pi (x+y)\sqrt{\frac nN})
\notag\\
&+\frac 1{2N}\sum_{n=1}^\infty\frac{1}{\lambda^{-1}e^{n/cN}+1}
\left(\frac nN\right)^{-1/2}\cos(\pi (x-y)\sqrt{\frac nN})
\to\frac 12\int_0^\infty\frac{\cos\pi(x-y)\sqrt{t}}{\lambda^{-1}e^{t/c}+1}\frac{dt}{\sqrt{t}}
\notag
\end{align}
uniformly for $x,y$ in a compact set as $N\to\infty$. 
If we make the change of variables 
$t=u^2$ we obtain $L_c(x,y)$ in (\ref{1.32}).
\end{proof}

When we are in the transition region in the bulk, the behaviour at the edge is still like 
that of independent random variables. This is the content of theorem \ref{thm1.9} which we 
prove next.

\begin{proof} (Theorem \ref{thm1.9}). 
We split the kernel $K_\lambda$ as follows
\begin{align}
K_\lambda(x,y)&=\sum_{n=0}^\infty A_N^{-1}e^{-n/cN}\psi_n(x)\psi_n(y)
\notag\\
&+\left(\sum_{n=0}^{M_N}+\sum_{n=M_N+1}^\infty\right)
\left(\frac 1{A_Ne^{n/cN}+1}-\frac 1{A_Ne^{n/cN}}\right)\psi_n(x)\psi_n(y)
\notag\\
&\doteq S_1(x,y)+S_2(x,y)+S_3(x,y),
\notag
\end{align}
where $M_N=[(1-\delta)cN\log N]$ with $\delta>0$ small. Here $A_N$ has the same
meaning as in the proof of theorem \ref{thm1.8}. Note that we have the estimate
\begin{equation}\label{4.5}
\left|\frac 1{A_Ne^{n/cN}+1}-\frac 1{A_Ne^{n/cN}}\right|\le Ce^{-2n/cN}.
\end{equation}
By (\ref{1.26}),
\begin{align}
&b_NS_1(a_N(c)+b_N\xi,a_N(c)+b_N\eta)
\notag\\
&=\frac{b_N}{A_N(1-q)\sqrt{\pi}}\exp(-\frac 12(a_N(c)+b_N\xi)^2-
\frac 12(a_N(c)+b_N\eta)^2-\frac{qb_N^2}{(1-q)^2}(\xi-\eta)^2)
\notag\\
&=\exp(-\frac{\xi+\eta}2-(c^2N^2+O(N))(\xi-\eta)^2 +o(1))
\notag
\end{align}
as $N\to\infty$. The last identity explains the choice of $a_N(c)$ and $b_N$ and we will
get (\ref{1.33}) if we can prove that $S_2$ and $S_3$ both tend to zero as $N$
tends to infinity.

From the estimates (\ref{4.3'}), (\ref{4.5}) and $\beta_q\sim\sqrt{2cN}$ we obtain
\begin{equation}
|b_NS_3(a_N(c)+b_N\xi,a_N(c)+b_N\eta)|\le \frac{CN^{1/2}}{\sqrt{\log N}}
\sum_{n=M_N+1}^\infty e^{-2n/cN}\frac 1{n^{1/6}}\to 0
\notag
\end{equation}
as $N\to\infty$ provided $\delta$ is sufficiently small.

We can write $\psi_n(a_N(c)+b_N\xi)=\psi_n(\sqrt{2n}y)$, where
\begin{equation}
y=\frac 1{\sqrt{2n}}(\sqrt{2cN}+O(\frac 1{N^{3/2}}))
(\sqrt{\log N}-(4\sqrt{\log N})^{-1}\log(\frac{4\pi}{\lambda^2c^2}\log N)
+\xi(2\sqrt{\log N})^{-1})
\notag
\end{equation}
For a fixed $\xi$,  we see that $y\ge 1+\delta$ if $N$ is sufficiently large and $1\le n\le M_N$. 
We can then use the estimate,
\begin{equation}\label{4.5'}
|h_n(\sqrt{2n}x)e^{-nx^2}|\le\frac{C_1}{n^{1/4}}e^{-nF(x)}
\end{equation}
for $x\ge 1+\delta$, \cite{DKMVZ}, which gives
\begin{equation}
|h_n(\sqrt{2n}x)e^{-nx^2}|\le\frac{C_1}{n^{1/4}}e^{-C_2n\delta^{3/2}},
\notag
\end{equation}
where $C_2$ is a numerical constant. If $N$ is large enough, then $y\ge\sqrt{cNn^{-1}\log N}$
and we get
\begin{align}
&|b_NS_2(a_N(c)+b_N\xi,a_N(c)+b_N\eta)|
\notag\\
&\le\frac 1{2\sqrt{N}}+\frac{C\sqrt{N}}{\sqrt{\log N}}\sum_{n=1}^{M_N}
\frac 1{n^{1/4}}e^{-C_2n(\sqrt{cNn^{-1}\log N}-1)^{3/2}}
\notag
\end{align}
which $\to 0$ as $N\to\infty$. 
\end{proof}

The next result shows that the kernel $M_\alpha$ can be obtained as a scaling limit of the
kernel $K_\lambda$

\begin{proof}
(Theorem \ref{thm1.10}).
Let $A_N$ have the same meaning as in the proof of theorem \ref{thm1.8}. We have
\begin{equation}
A_N\approx e^{-\alpha N^{2/3}+\alpha N^{-1/3}/2}
\notag
\end{equation}
with a negligible error. Also, as $N\to\infty$,
\begin{equation}
\beta_q=\frac{\sqrt{2}N^{1/6}}{\sqrt{\alpha}}+O(\frac 1{\sqrt{N}}).
\notag
\end{equation}
Write $f_n(x)=h_n(x)e^{-x^2/2}$ as above. We have
\begin{align}
&\frac{\sqrt{\alpha}}{2N^{1/3}}K_\lambda(N^{1/3}\sqrt{\alpha}+
\frac{\sqrt{\alpha}}{2N^{1/3}}\xi,N^{1/3}\sqrt{\alpha}+
\frac{\sqrt{\alpha}}{2N^{1/3}}\eta)
\notag\\
&=\left(\frac 1{N^{1/6}\sqrt{2}}+O(\frac 1{N^{5/6}})\right)
\sum_{n=0}^\infty\frac 1{e^{\alpha N^{2/3}((n+1/2)/N-1)}+1}
\notag\\
&\times f_n(\sqrt{2N}+\frac{\xi}{N^{1/6}\sqrt{2}}+O(\frac 1{N^{1/3}}))
f_n(\sqrt{2N}+\frac{\eta}{N^{1/6}\sqrt{2}}+O(\frac 1{N^{1/3}}))
\notag\\
&=\left(\frac 1{N^{1/6}\sqrt{2}}+O(\frac 1{N^{5/6}})\right)
\left(\sum_{k=-\infty}^{-M_N-1}+\sum_{k=-M_N}^{M_N}+\sum_{k=M_N+1}^N\right)
\frac 1{e^{-\alpha(k-1/2)/N^{1/3}}+1}
\notag\\
&\times f_{N-k}(\sqrt{2N}+\frac{\xi}{N^{1/6}\sqrt{2}}+O(\frac 1{N^{1/3}}))
f_{N-k}(\sqrt{2N}+\frac{\eta}{N^{1/6}\sqrt{2}}+O(\frac 1{N^{1/3}}))
\notag\\
&\doteq \Sigma_1+\Sigma_2+\Sigma_3,
\notag
\end{align}
where $M_N=\gamma N^{1/3}\log N$ with some fixed $\gamma>0$ that can be chosen.
The asymptotic contribution will come from $\Sigma_2$. Here we use the asymptotic formula
\begin{equation}\label{4.6}
f_{N-k}(\sqrt{2N}+\frac{u}{N^{1/6}\sqrt{2}})=\frac{2^{1/4}}{N^{1/12}}
\Ai(u+\frac{k-1/2}{N^{1/3}})(1+O(\frac{\log N}{N^{2/3}}))
\end{equation}
for $|k|\le\gamma N^{1/3}\log N$ and $u$ in a compact set. This formula folows from results
in \cite{DKMVZ}, see \cite{AdvM}. Using this we see that
\begin{align}
\lim_{N\to\infty} \Sigma_2&=\lim_{N\to\infty}\frac 1{N^{1/3}}\sum_{-M_N}^{M_N}
\frac 1{e^{-\alpha(k-1/2)/N^{1/3}}+1}\Ai(\xi+\frac{k-1/2}{N^{1/3}})
\Ai(\eta+\frac{k-1/2}{N^{1/3}})
\notag\\
&=\int_{-\infty}^\infty\frac 1{e^{-\alpha x}+1}\Ai(x+\xi)\Ai(x+\eta)dx=M_\alpha(\xi,\eta).
\notag
\end{align}

We still have to prove $\Sigma_1\to 0$ and $\Sigma_3\to 0$ as $N\to\infty$.
For $\Sigma_1$ we use the estimate (\ref{4.3'}), which gives
\begin{equation}
|\Sigma_1|\le\frac C{N^{1/6}}\sum_{k=-\infty}^{-M_N-1}
\frac 1{e^{-\alpha(k-1/2)/N^{1/3}}+1}\frac 1{(N-k)^{1/6}}
\notag
\end{equation}
which goes to zero as $N$ tends to infinity if we choose $\gamma$ large enough.

For $1<x\le 1+\delta$ we have an asymptotic formula for $f_n(\sqrt{2n}x)$ in terms of
the Airy function, see \cite{DKMVZ}. Estimates of the Airy function then gives
\begin{equation}\label{4.7}
|f_n(\sqrt{2n}x)|\le\frac C{N^{1/12}}e^{-cn(x-1)^{3/2}}
\end{equation}
for some constants $c,C$, when $1<x\le 1+\delta$. Since $F(x)\ge c(x-1)^{3/2}$, we can
combine this with (\ref{4.5'}) to see that (\ref{4.7}) holds for all $x>1$. For $N$ 
sufficiently large this leads to an estimate
\begin{equation}
\left|f_{N-k}(\sqrt{2N}+\frac{\xi}{N^{1/6}\sqrt{2}}+O(\frac 1{N^{1/3}}))\right|
\le Ce^{-ck^{3/2}/\sqrt{N}}.
\notag
\end{equation}
It follows that 
\begin{equation}
|\Sigma_3|\le C\sum_{k=M_N+1}^N e^{-ck^{3/2}/\sqrt{N}}\to 0
\notag
\end{equation}
as $N\to\infty$ if we choose $\gamma$ sufficiently large.
\end{proof}

We give here a sketch of an argument for the approximate expression (\ref{1.32'}) for 
$L_c(x,y)$. Integration by parts gives
\begin{align}
L_c(x,y)&=\int_0^\infty\frac{\sin\pi(x-y)u}{\pi (x-y)}\frac {2u}{\lambda c}
\frac{e^{u^2/c}}{(\lambda^{-1}e^{u^2/c}+1)^2}du
\notag\\
&=\frac 2c\int_0^\infty u\frac{\sin\pi(x-y)u}{\pi (x-y)}\frac{du}{\cosh^2
(u^2/2c+a/2)},
\notag
\end{align}
where $a=\log(1/\lambda)$. When $c$ is small $a\approx -1/c$.
Make the change of variables $u=1+ct$. This gives
\begin{align}
L_c(x,y)&=\int_{-1/c}^\infty(1+ct)\frac{\sin\pi(x-y)(1+ct)}{\pi (x-y)}
\frac{dt}{\cosh^2(t+ct^2/2)}
\notag\\
&\approx\frac 12\int_{-\infty}^\infty (1+ct)\frac{\sin\pi(x-y)(1+ct)}{\pi (x-y)}
\frac{dt}{\cosh^2 t}.
\notag
\end{align}
If we use the addition formula for the sine function and neglect terms containg $c^2$
we get
\begin{align}
L_c(x,y)&\approx\int_{-\infty}^\infty\frac{\cos\pi(x-y)ct}{\cosh 2t+1}
\frac{\sin\pi(x-y)}{\pi (x-y)}
\notag\\
&=\frac{\pi^2(x-y)c}{2\sinh(\pi^2c(x-y)/2)}\frac{\sin\pi(x-y)}{\pi (x-y)}=
\frac{\pi c}2\frac{\sin\pi(x-y)}{\sinh(\pi^2c(x-y)/2)}.
\notag
\end{align}
We see that as $c\to 0+$ the kernel $L_c$ approaches the sine kernel.

\section{Largest eigenvalue for deformed GUE}\label{sect5}

This section contains the proof of theorem \ref{thm1.11}.
Consider $N$ non-intersecting standard Brownian motions started at $y_1,\dots,y_N$ and
conditioned never to intersect. If we fix $y_1,\dots,y_N$ the particle distribution at
time $S$ is a determinantal process with correlation functions
\begin{equation}\label{5.1}
\rho_{m,N}(x_1,\dots,x_m;y)=\det(K_{N}(x_i,x_j;y))_{1\le i,j\le m},
\end{equation}
where
\begin{equation}\label{5.2}
K_N(u,v;y)=\frac 1{(2\pi i)^2S}\int_\gamma dz\int_\Gamma e^{(w-v)^2/2S-(z-u)^2/2S}
\frac 1{w-z}\prod_{j=1}^N\frac{w-y_j}{z-y_j},
\end{equation}
see e.g. \cite{JoUni}. Here $\gamma$ is a positively oriented simple closed curve
containing $y_1,\dots,y_N$ and $\Gamma$ a verical line oriented upwards and not interesecting 
$\gamma$, we place it to the right of $\gamma$. The formula (\ref{5.1}) also gives the
correlation functions for the eigenvalues of the hermitian matrix $M$ given by
\begin{equation}\label{5.3}
M=\text{diag\,}(y_1+\dots,y_N)+\sqrt{2S} V
\end{equation}
with $V$ a standard $N\times N$ GUE matrix. If $\lambda_{\max}^{(N)}$ is the largest 
eigenvalue of $M$ then $\lambda_{\max}^{(N)}$ has the same distribution as
$\max_{1\le j\le N} x_j$.

Let $\mathbb{P}_y$ denote the probability measure for $y_1,y_2,\dots$ and let $\mathbb{P}_
{x;y}$ denote the expectation with respect to the determinantal process with correlation
kernel $K_N(u,v;y)$ given by (\ref{5.2}). Furthermore, we 
let $\mathbb{P}_N=\mathbb{P}_y\otimes\mathbb{P}_{x;y}$ be the product measure.
We are interested in the distribution function
\begin{align}\label{5.4}
F_N(t)&=\mathbb{P}[\lambda_{\max}^{(N)}\le t]=\mathbb{E}_N[\prod_{j=1}^N(1-\chi_{(t,\infty)}
(x_j))]
\notag\\
&=\mathbb{E}_y[\mathbb{E}_{x;y}[\prod_{j=1}^N(1-\chi_{(t,\infty)}(x_j))]].
\end{align}
When computing the inner expectation we are considering
$y_1,\dots, y_N$ as fixed and hence we can work with the correlation functions (\ref{5.1}). 

Fix a number $\epsilon\in(1/7,1/6)$ and set
\begin{equation}
A_N=\{y\in\mathbb{R}^N\,;\,|y_i|\le N^\epsilon\,,\,1\le i\le N\}.
\notag
\end{equation}
Define a cut-off measure $d\mu_N(t)$ by
\begin{equation}\label{5.5}
d\mu_N(t)=\frac 1{\mu([-N^\epsilon,N^\epsilon])}\chi_{[-N^\epsilon,N^\epsilon]}(t)d\mu(t),
\end{equation}
and the function
\begin{equation}\label{5.6}
G_N(z)=\int_{\mathbb{R}}\frac{d\mu_N(t)}{z-t}
\end{equation}
for $z\in\mathbb{C}\setminus [-N^\epsilon,N^\epsilon]$.

Since by assumption $\mu$ has finite 7:th moment it follows that
\begin{equation}
\mathbb{P}_y[A_N^c]\le N\frac{C}{N^{7\epsilon}}
\notag
\end{equation}
which $\to 0$ as $N\to\infty$. Hence, since the expression in the $\mathbb{E}_y$-expectation 
in (\ref{5.4}) is bounded, we can restrict our attention to $A_N$ and use $\mu_N$ instead
of $\mu$, so we regard $y_1,\dots,y_N$ as independent random variables with distribution
$\mu_N$. Denote this probability measure by $\mathbb{P}_y^{(N)}$.

\begin{lemma}\label{lem5.1}
There is a real number $w_c=w_c(n)$, which is approximately $\sqrt{NS}=\alpha N^{1/6}$,
such that
\begin{equation}\label{5.7}
G_N'(w_c)=-\frac 1{\alpha^2N^{1/3}}
\end{equation}
for all sufficiently large $N$.
\end{lemma}

\begin{proof}
We have
\begin{equation}
G_N'(z)=-\int_{-N^\epsilon}^{N^\epsilon}\frac{d\mu_N(t)}{(z-t)^2}.
\notag
\end{equation}
The moment conditions on $\mu$ can be used to see that
\begin{equation}\label{5.8}
|G_N'(z)-\frac 1{z^2}|\le\frac C{z^3}
\end{equation}
for real $z\ge 2N^\epsilon$ say. We see that, for $N$ sufficiently large, $G_N'$ is
a decreasing function in $[2N^\epsilon,\infty)$. Furthermore, $G_N'(z)\to 0$ as $z\to\infty$
and $G_N''(2N^\epsilon)\approx 1/4N^{2\epsilon}>1/\alpha N^{1/3}$
Hence, when $N$ is large enough, there is a $z=w_c$ such that (\ref{5.7}) holds.
From (\ref{5.8}) we see that $1/w_c^2\sim1/\alpha^2N^{1/3}$, which gives
the asymptotic behaviour.
\end{proof}

Set 
\begin{equation}\label{5.9}
r_N(y)=-\sum_{j=1}^N\frac 1{(w_c-y_j)^2}+NG_N'(w_c),
\end{equation}
and
\begin{equation}\label{5.10}
v_c=w_c+S\sum_{j=1}^N\frac 1{w_c-y_j}.
\end{equation}
Furthermore, set
\begin{equation}\label{5.11}
R(N)=w_c+\frac{\alpha^2N^{1/3}}{w_c}+\frac{\alpha^2N^{1/3}}{w_c}\int_{\mathbb{R}}
\frac y{1-y/w_c}d\mu_N(y)
\end{equation}
and define
\begin{equation}\label{5.12}
s_N(y)=\frac{\alpha}{w_cN^{1/6}}\left(\sum_{j=1}^N\frac {y_j}{w_c-y_j}-
\int_{\mathbb{R}}\frac y{w_c-y}d\mu_N(y)\right).
\end{equation}
Note that
\begin{equation}\label{5.13}
v_c=R(N)+\frac{\alpha}{N^{1/2}}s_N(y).
\end{equation}

\begin{lemma}\label{lem5.2}
\item{(i)} There is a constant $C$ such that
\begin{equation}
\mathbb{E}_y^{(N)}[r_N(y)]\le C.
\end{equation}
\item{(ii)} We have the limit
\begin{equation}
\text{Var}_y^{(N)}[s_N(y)]\to\sigma^2/\alpha^2
\end{equation}
as $N\to\infty$.
\item{(iii)} The random variable $s_N(y)$ converges in distribution to
$N(0,\sigma^2/\alpha^2)$ as $N\to\infty$
\end{lemma}

\begin{proof} Since the 7:th moment is finite and $\int td\mu(t)=0$ 
we get $|\int td\mu_N(t)|\le CN^{-6\epsilon} \le CN^{-6/7}$. The definition of $w_c$ and 
$G_N$ gives $\mathbb{E}_y^{(N)}[r_N(y)]=0$ and thus
\begin{align}
\mathbb{E}_y^{(N)}[r_N(y)^2]&=\text{Var}_y^{(N)}[r_N(y)]=
\text{Var}_y^{(N)}[\sum_{j=1}^N\frac 1{(w_c-y_j)^2}]
\notag\\
&=N[-\frac 16G_N^{(3)}(w_c)-G_N'(w_c)^2]\le\frac {NC}{w_c^6},
\notag
\end{align}
where the last inequality follows from our moment condition. Since $w_c\sim\alpha N^{1/6}$
we see that the right hand side is bounded. To prove (ii) we compute
\begin{equation}\label{5.13'}
\text{Var}_y^{(N)}[s_N(y)]=\frac{\alpha^2 N^{2/3}}{w_c^4}\left(
\int\frac{y^2}{(1-y/w_c)^2}d\mu_N(y)-\left(\int\frac{y}{1-y/w_c}d\mu_N(y)\right)^2\right).
\end{equation}
Writing $y(1-y/w_c)^{-1}=y+y^2/w_c+\dots$, we get $\int y(1-y/w_c)^{-1}d\mu_N(y)=o(1)$. 
Similarly, writing $y^2/(1-y/w_c)^2=y^2+2y^3/w_c+\dots$, we get
$\int y^2/(1-y/w_c)^2d\mu_N(y)=\sigma^2+o(1)$.
Using $w_c\sim\alpha N^{1/6}$, the identity (\ref{5.13'}) now yields (ii). The claim
(iii) follows immediately from the central limit theorem, since $s_N$ is a sum of 
independent random variables.
\end{proof}

Define $B_N$ to be the set of all $y\in A_N$ such that $|r_N(y)|\le CN^\epsilon$ and
$|s_N(y)|\le CN^\epsilon$. It follows from lemma \ref{lem5.2} (i) and (ii) that we can
restrict ourselves to $y\in B_N$. We will study
\begin{equation}\label{5.14}
F_N^\ast(t)=\mathbb{E}_y^{(N)}[\chi_{B_N}(y)\mathbb{E}_{x;y}[
\prod_{j=1}^N(1-\chi_{(t,\infty)}(x_j))]]
\end{equation}
instead of $F_N(t)$ given by (\ref{5.4}). Hence it is enough to consider a fixed
$y$ in $B_N$ and work with
\begin{align}\label{5.15}
F_N^\ast(t;y)&=\mathbb{E}_{x;y}[\prod_{j=1}^N(1-\chi_{(t,\infty)}(x_j))]
\notag\\
&=\sum_{k=0}^N\frac{(-1)^k}{k!}\int_{(t,\infty)^k}\det(K_N(x_i,x_j;y))d^kx.
\end{align}

Our problem is then to investigate the asymptotics of $K_N(u,v;y)$ for $y\in B_N$.
Set
\begin{equation}
f(w)=\frac{w^2}{2S}-\frac {v_cw}S+\sum_{j=1}^N\log(w-y_j),
\notag
\end{equation}
where we choose the principal branch of the logarithm. The number $v_c$
was defined so that
\begin{equation}\label{5.16}
f'(w_c)=0.
\end{equation}
We consider
\begin{equation}\label{5.17}
u=v_c+\xi\frac{\alpha}{\sqrt{N}}\,\,;\,\,v=v_c+\eta\frac{\alpha}{\sqrt{N}},
\end{equation}
where $\xi$ and $\eta$ lie in a compact set.

To perform a saddle-point argument in the integral (\ref{5.2}) we must specify 
appropriate contours. Let $C_1:[0,\infty)\ni t\to w_c+t+it$,
$C_2:[0,\infty)\ni t\to w_c+t-it$, $C_3:[0,\infty)\ni t\to w_c-t+it$
and $C_4:[0,\infty)\ni t\to w_c-t-it$. We want to show that we can deform 
$\Gamma$ to $C_1-C_2$ and $\gamma$ to $C_3-C_4$ in the contour integral
(\ref{5.2}).
Let $C^A_i$ be the parts of the contours where we restrict $t$ to $[0,A+w_c]$
and let $\gamma_A:[-w_c-A,w_c+A]\ni t \to -A-it$, where $A>0$. Then
$\gamma$ can be deformed to $C^A_3+\gamma_A-C^A_4$ in (\ref{5.2}) if $A$ is
sufficiently large. From (\ref{5.11}) we see that $R(N)\sim 2\alpha N^{1/6}$ and since
$\xi$ belongs to a compact set we see from (\ref{5.13}), $|s_N(y)|\le N^\epsilon$ and
(\ref{5.17}) that $u\ge \alpha N^{1/6}$ for all sufficiently large $N$. Since
$\re (-(z-u)^2)\le -2uA$ we see that we can let $A\to\infty$ and conclude that the
contribution from $\gamma_A$ goes to zero.

Choose $\Gamma$ to be a vertical line through $w_c$. We want to show that 
the part of this line that lies in the upper half plane can be deformed to
$C_1$, and the in the lower half plane to $C_2$. Set $w=w_c+t+iA$, $0\le t\le A$.
Then,
\begin{equation}
g(t)=\re((w-v)^2)=t^2+2(w_c-v)t+(w_c-v)^2-A^2.
\notag
\end{equation}
For $0\le t\le A/2$, we see that $g(t)\le -A^2/2$ for $A$ large, 
and when $A/2\le t\le A$ we have
\begin{equation}
g(t)\le (w_c-v)A+(w_c-v)^2\le -\frac{\alpha}2 N^{1/6}A
\notag
\end{equation}
for $N$ and $A$ large. Hence we can deform the upper part of $\Gamma$ to $C_1$ in 
(\ref{5.2}). The deformation to $C_2$ is analogous by symmetry.

Next, we want to localize the integration to a small neighbourhood of $w_c$.

\begin{lemma}\label{lem5.3}
Define $g_i(t)=\re(f(C_i(t))-f(w_c))$ for $i=1,2$ and \linebreak
$g_i(t)=-\re(f(C_i(t))-f(w_c))$ for $i=3,4$, $t\ge 0$. There is a positive constant
$c$ so that
\begin{equation}\label{5.18}
g_i(t)\le\begin{cases} -cN^{1/2}t^3 &,0\le t\le\alpha N^{1/6}/2\\
-cN^{5/6}t & t\ge\alpha N^{1/6}/2\end{cases},
\end{equation}
for all sufficiently large $N$.
\end{lemma}

\begin{proof}
Consider $g_1(t)$. We have
\begin{equation}
g_1(t)=\frac 1{2S}[(w_c+t)^2-t^2]-\frac 12 v_c(w_c+t)+
\frac 12\sum_{j=1}^N\log((w_c+t-y_j)^2+t^2)-\re f(w_c).
\notag
\end{equation}
Differentiation gives
\begin{equation}
g_1'(t)=\frac{w_c-v_c}S+\sum_{j=1}^N\frac{w_c-y_j+2t}{(w_c-y_j+t)^2+t^2}
\notag
\end{equation}
and hence $g_1'(0)=0$ by (\ref{5.10}). Hence
\begin{equation}
g_1'(t)=g_1'(t)-g_1'(0)=-\sum_{j=1}^N\frac{2t^2}{(w_c-y_j)((w_c-y_j+t)^2+t^2)}.
\notag
\end{equation}
We know that $w_c\sim\alpha N^{1/6}$ and $|y_j|\le N^\epsilon$. If 
$|t|\le\alpha N^{1/6}/2$, we see that there is a positive constant $c$ such that 
$g_1'(t)\le-cN^{1/2}t^2$ for $0\le t\le\alpha N^{1/6}/2$. Since $g_1(0)=0$ we obtain the 
first part of (\ref{5.18}) for $i=1$. If $t\ge \alpha N^{1/6}/2$, then $g_1'(t)\ge
-cN^{5/6}$ for some postive constant $c$, and we obtain the second part of (\ref{5.18}).

Consider next $g_3(t)$. Again $g_3(0)=g_3'(0)=0$ and we get
\begin{equation}
g_3'(t)=-\sum_{j=1}^N\frac{2t^2}{(w_c-y_j)((w_c-y_j+t)^2+t^2)}
\notag
\end{equation}
and we can proceed as above. The functions $g_2$ and $g_3$ are treated analogously.
\end{proof}

We also need a local approximation of $f(w)$ in a neighbourhood of $w_c$. By (\ref{5.16})
we have $f'(w_c)=0$ and we also have

\begin{equation}\label{5.19}
f''(w_c)=\frac 1S-\sum_{j=1}^N\frac 1{(w_c-y_j)^2}=r_N(y),
\end{equation}
by (\ref{5.7}) and (\ref{5.9}). Furthermore,
\begin{equation}\label{5.20}
f^{(3)}(w_c)=\sum_{j=1}^N\frac 2{(w_c-y_j)^3}.
\end{equation}

\begin{lemma}\label{lem5.4}
For $\zeta\in\mathbb{C}$ and $|\zeta|\le N^{1/18}$,
\begin{equation}\label{5.21}
f(w_c+\zeta\alpha N^{-1/6})=f(w_c)+\frac 13\zeta^3+o(1),
\end{equation}
as $N\to\infty$, where $o(1)$ is uniform for $|\zeta|\le N^{1/18}$.
\end{lemma}

\begin{proof}
Define $R(\zeta)$ by
\begin{equation}\label{5.22}
f(w_c+\zeta)=f(w_c)+f'(w_c)\zeta+\frac 12f''(w_c)\zeta^2+\frac 16 f^{(3)}(w_c)\zeta^3
+R(\zeta).
\end{equation}
Since $w_c\sim\alpha N^{1/6}$ and $|w_c-y_j|$ is much greater than 1 for $N$ large enough 
a Taylor expansion gives
\begin{equation}\label{5.23}
|R(\zeta)|\le CN^{1/3}|\zeta|^4
\end{equation}
for $|\zeta|\le 1$. It follows from (\ref{5.16}), (\ref{5.19}) and (\ref{5.20}) that
\begin{equation}\label{5.24}
f(w_c+i\zeta\alpha N^{-1/6})=f(w_c)-\frac{\alpha^2}2 r_N(y)\zeta^2N^{-1/3}
+\frac 13\sum_{j=1}^N\frac {\zeta^3N^{-1/2}} {(w_c-y_j)^3}+R(\zeta\alpha N^{-1/6}).
\end{equation}
Since $|r_N(y)|\le CN^\epsilon$ and $|\zeta|\le N^{1/18}$ we see that
$|r_N(y)\zeta^2N^{-1/3}|\le CN^{-1/18}$. Furthermore, by (\ref{5.23}), 
$|R(\zeta\alpha N^{-1/6})|\le CN^{-1/9}$. We can write
\begin{equation}
\sum_{j=1}^N\frac 1{(w_c-y_j)^3}=\frac N{w_c^3}+\sum_{j=1}^N
\frac{w_c^3-(w_c-y_j)^3}{w_c^3(w_c-y_j)^3}=\frac{N^{1/2}}{\alpha^3}(1+o(1))
\notag
\end{equation}
as $N\to\infty$. We see now that (\ref{5.21}) follows from (\ref{5.24}).
\end{proof}

It follows from (\ref{5.2}), the definition of $f$ and the change of contours discussed
above that
\begin{align}\label{5.25}
&K_N(v_c+\xi\frac{\alpha}{\sqrt{N}},v_c+\eta\frac{\alpha}{\sqrt{N}})
\notag\\
&=\frac
{e^{v^2-u^2}N^{2/3}}{(2\pi i)^2\alpha^2}\int_{C_3-C_4}dz\int_{C_1-C_2}dw
\frac{e^{f(w)-f(z)}}{w-z}e^{-\eta N^{1/6}w/\alpha+\xi N^{1/6} z/\alpha}
\end{align}
Consider $z$ on $C_3$ and $w$ on $C_1$. The other cases are similar. Set
$z=w_c+(-t+it)\alpha N^{-1/6}$, $w=w_c+(-\tau+i\tau)\alpha N^{-1/6}$,
$t,\tau\ge 0$. It follows from lemma \ref{lem5.3} that we can localize the
evaluation of (\ref{5.25}) to $t,\tau\le N^{1/18}$. By lemma \ref{lem5.4}
\begin{equation}
f(w)-f(z)=\frac 13(\tau+i\tau)^3-\frac 13(-t+it)^3+o(1)
\notag
\end{equation}
uniformly for $0\le t,\tau\le N^{1/18}$. Hence, the contribution to (\ref{5.25}) from
$z$ on $C_3$ and $w$ on $C_1$ is
\begin{equation}
\frac{N^{1/2}e^{v^2-u^2+(\xi-\eta)N^{1/6}w_c/\alpha}}{(2\pi i)^2\alpha}
\int_0^{N^{1/18}}dt\int_0^{N^{1/18}}d\tau
\frac{e^{\frac 13(\tau+i\tau)^3-\frac 13(-t+it)^3+\xi(-t+it)-\eta(\tau+i\tau)}} 
{(\tau+i\tau)-(-t+it)}.
\notag
\end{equation}
Define
\begin{equation}\label{5.26}
K_N^\ast(u,v;y)= e^{v^2-u^2+(u-v)w_c/S}K_N(u,v;y).
\end{equation}
We can just as well use $K_N^\ast$ as $K_N$. If we argue as above for all parts of
the contours we get
\begin{align}\label{5.27}
&\lim_{N\to\infty}\frac{\alpha}{\sqrt{N}}K_N^\ast(v_c+\xi\frac{\alpha}{\sqrt{N}},
v_c+\eta\frac{\alpha}{\sqrt{N}})
\notag\\
&=\frac 1{(2\pi i)^2}\int_{\gamma'}dz\int_{\Gamma'}dw e^{w^3/3-z^3/3+\xi z-\eta w}\frac
1{w-z},
\end{align}
where $\gamma'$ is the contour given by $t+it$ for $t\le 0$ and $-t+it$ for $t\ge 0$,
and $\Gamma'$ is the reflection of $\gamma'$ in the imaginary axis.

We have to show that the right hand side of (\ref{5.27}) is really the Airy kernel.
Let $\mathcal{C}$ be the contour given by $t+i|t|$, $t\in\mathbb{R}$. 
If we change variables
by $z=i\zeta$, $w=-i\omega$, then $\gamma'$ maps to $\mathcal{C}$ and $\Gamma'$ to
$-\mathcal{C}$ and we see that the right hand side of (\ref{5.27}) becomes
\begin{equation}
-\frac 1{4\pi^2}\int_{\mathcal{C}}d\zeta\int_{\mathcal{C}}d\omega e^{i\omega^3/3+i\eta 
\omega+i\zeta^3/3+i\xi\zeta}\frac 1{i(\zeta+\omega)},
\notag
\end{equation}
which is the Airy kernel $K_{Airy}(\xi,\eta)$, \cite{JoDPNG}.

We have proved

\begin{lemma}\label{lem5.5}
Define
\begin{equation}\label{5.27'}
\tilde{K}_N(\xi,\eta;y)=\frac{\alpha}{\sqrt{N}}K_N^\ast(v_c+\xi\frac{\alpha}{\sqrt{N}},
v_c+\eta\frac{\alpha}{\sqrt{N}})
\end{equation}
with $K_N^\ast$ as in (\ref{5.26}) and $v_c$ given by (\ref{5.13}). Then
\begin{equation}\label{5.28}
\lim_{N\to\infty}\tilde{K}_N(\xi,\eta;y)=K_{Airy}(\xi,\eta)
\end{equation}
uniformly for $\xi,\eta$ in a compact set and $y\in B_N$.
\end{lemma}

To control the convergence of (\ref{5.15}) we need some more estimates.

\begin{lemma}\label{lem5.6}
Fix a constant $A$. There is a constant $C$, depending on $A$, such that for
$\xi,\eta\ge -A$ and all sufficiently large $N$ we have the estimate
\begin{equation}\label{5.29}
\tilde{K}_N(\xi,\eta;y)|\le Ce^{-\frac 23(|\xi|^{3/2}+|\eta|^{3/2})}
\end{equation}
for $y\in B_N$.
\end{lemma}

\begin{proof}
Deform the contour $C_3-C_4$ to $\gamma=\gamma_1+\gamma_2+\gamma_3$ and
$C_1-C_2$ to $\Gamma=\Gamma_1+\Gamma_2+\Gamma_3$, where
$\gamma_1:(-\infty,-\alpha\delta N^{-1/6})\ni t\to w_c+t+it$,
$\gamma_2:(-\delta,\delta)\ni t\to w_c-\alpha\delta N^{-1/6}+\alpha itN^{-1/6}$,
$\gamma_3:(\alpha\delta N^{-1/6},\infty)\ni t\to w_c-t+it$ and
$\Gamma_1,\Gamma_2,\Gamma_3$ are obtained by reflection in the line $\re z=w_c$.

Note that for $z$ on $\gamma$ and $w$ on $\Gamma$,
\begin{equation}
\frac 1{|w-z|}\le\frac{N^{1/6}}{2\alpha\delta}
\notag
\end{equation}
From (\ref{5.25}), (\ref{5.26}) and (\ref{5.27'}) we obtain
\begin{align}\label{5.29'}
\tilde{K}_N(\xi,\eta;y)|&\le \frac{N^{1/3}}{8\alpha^2\delta\pi^2}
\left(\int_\gamma e^{-\re(f(z)-f(w_c))+\xi N^{1/6}\alpha^{-1}\re(z-w_c)}|dz|\right)
\notag\\
&\times\left(\int_\gamma e^{-\re(f(w)-f(w_c))+\eta N^{1/6}\alpha^{-1}\re(w-w_c)}
|dw|\right).
\end{align}
On $\Gamma_2$ we can use lemma \ref{lem5.4} to get
\begin{equation}
\re[f(w_c+\alpha(\delta+it)N^{-1/6})-f(w_c)]=
\frac 13\delta^3-\delta t^2+o(1).
\notag
\end{equation}
This gives
\begin{align}
&\int_{\Gamma_2} e^{-\re(f(w)-f(w_c))+\eta N^{1/6}\alpha^{-1}\re(w-w_c)}
|dw|
\notag\\
&\le\frac{C}{N^{1/6}}e^{\delta^3/3-\eta\delta}\int_{-\delta}^\delta e^{-\delta t^2}dt
\le\frac{C}{N^{1/6}\delta^{1/2}}e^{\delta^3/3-\eta\delta}.
\notag
\end{align}
On $\Gamma_3$ we can use lemma \ref{lem5.3} to get
\begin{align}
&\int_{\Gamma_3} e^{-\re(f(w)-f(w_c))+\eta N^{1/6}\alpha^{-1}\re(w-w_c)}
|dw|
\notag\\
&\le\int_{-\delta\alpha N^{-1/6}}^{\alpha N^{1/6}/2} 
e^{-cN^{1/2}t^3-\eta N^{1/6}t/\alpha}dt+
+\int_{\alpha N^{1/6}/2}^\infty e^{-c N^{5/6} t-\eta N^{1/6}t/\alpha}dt
\notag\\
&\le\frac C{N^{1/6}}e^{-\eta\delta}.
\notag
\end{align}
The contribution from $\Gamma_1$ is analogous. Choosing $\delta=\sqrt{\eta}$ for $\eta
\ge 1$ and $\delta=1$ otherwise, we get
\begin{equation}
\int_\Gamma e^{-\re(f(w)-f(w_c))+\eta N^{1/6}\alpha^{-1}\re(w-w_c)}|dw|
\le\frac{C}{N^{1/6}}e^{-\frac 23|\eta|^{3/2}}.
\notag
\end{equation}
The estimate for the other integral in (\ref{5.29'}) is analogous and the estimate 
(\ref{5.29}) follows.
\end{proof}

Define the distribution function $H_N$ by
\begin{equation}\label{5.30}
H_N(t;y)=\sum_{k=0}^N\frac{(-1)^k}{k!}\int_{(t,\infty)^k}\det(\tilde{K}_N
(\xi_i,\xi_j;y))d^k\xi.
\end{equation}
It follows from lemma \ref{lem5.5} and lemma \ref{lem5.6} that
\begin{equation}\label{5.31}
\lim_{N\to\infty}H_N(t;y)=F_{TW}(t)
\end{equation}
uniformly for $t$ in a compact subset and $y\in B_N$.

If we change from $K_N$ to $K_N^\ast$ in (\ref{5.15}) and make the change of variables
$x_i=v_c+\alpha\xi_i/\sqrt{N}$, we see that
\begin{equation}
F_N^\ast(t)=\mathbb{E}_y^{(N)}\left[\chi_{B_N}(y)G_N(\sqrt{N}(t-v_c)/\alpha;y)\right].
\notag
\end{equation}
Thus
\begin{equation}
F_N^\ast(R(N)+\xi\alpha/\sqrt{N})=
\mathbb{E}_y^{(N)}[\chi_{B_N}(y)G_N(\xi-s_N(y);y)].
\notag
\end{equation}
We can now use lemma \ref{lem5.2} and (\ref{5.31}) to see that
\begin{equation}
\lim_{N\to\infty}F_N^\ast(R(N)+\xi\alpha/\sqrt{N})=
\int_{\mathbb{R}}F_{TW}(\xi-u)h(u)du,
\notag
\end{equation}
where $h(u)=(2\pi\sigma^2/\alpha^2)^{-1/2}\exp(-u^2\alpha^2/2\sigma^2)$. This 
completes the proof of theorem \ref{thm1.11}.

\medskip
\noindent
{\bf Acknowledgement}: I thank Y. Chen for drawing my attention to some papers
on intermediate ensembles.

\end{document}